\documentclass[11pt]{amsart}
\usepackage{amsfonts}
\usepackage{amsmath}
\usepackage{amssymb}
\usepackage{graphicx}
\usepackage{amsthm,amscd}
\usepackage{calc}
\usepackage{hyperref}

\newtheorem{theorem}{Theorem}[section]

\newtheorem{corollary}{Corollary}[section]

\newtheorem{definition}{Definition}[section]

\newtheorem{lemma}{Lemma}[section]

\newtheorem{proposition}{Proposition}[section]
\newtheorem{remark}{Remark}[section]

\numberwithin{equation}{section} \theoremstyle{remark}

\renewcommand{\bar}{\overline}

\newcommand{\pa}{\partial}
\renewcommand{\phi}{\varphi}

\newcommand{\ka}{K\"ahler }

\newcommand{\C}{{\mathbb C}}

\newcommand{\R}{{\mathbb R}}

\newcommand{\pz}{\partial_z}
\newcommand{\pzb}{\partial_{\bar z}}
\newcommand{\M}{{\mathcal M}}
\newcommand{\T}{{\mathcal T}}

\newcommand{\ke}{K\"ahler-Einstein }

\newcommand{\tei}{Teichm\"uller }
\newcommand{\h}{{\mathbb H}}
\renewcommand{\tilde}{\widetilde}
\newcommand{\p}{{\Phi}}

\newcommand{\lb}{\left (}
\newcommand{\rb}{\right )}
\newcommand{\po}{Poincar\'e }
\newcommand{\pg}{Poincar\'e growth }
\newcommand{\nspg}{Poincar\'e growth}

\title[]{Recent Development on the Geometry of the Teichm\"uller and Moduli Spaces of Riemann Surfaces and Polarized Calabi-Yau Manifolds}

\author{Kefeng Liu, Xiaofeng Sun, Shing-Tung Yau}

\begin{document}

\maketitle

\tableofcontents

\section{Introduction}\label{intro}

The moduli space $\M_{g,k}$ of Riemann surfaces of genus $g$ with
$k$ punctures plays an important role in many area of mathematics
and theoretical physics.
 In this article we first survey some of our recent works on the geometry of this moduli space.
  In the following we assume $g\geq 2$ and $k=0$ to simplify notations.
  All the results in this paper work for the general case when
  $3g-3+k>0$. We will focus on the K\"ahler metrics on the moduli
  and Teichm\"uller spaces, especially the Weil-Petersson metric,
  the Ricci, the perturbed Ricci, and the K\"ahler-Einstein metrics.

  We will review certain new geometric properties we found and proved for
  these metrics, such as the bounded geometry, the goodness and
  their naturalness under restriction to boundary divisors. The algebro-geometric corollaries such as the stability of the
  logarithmic cotangent bundles and the infinitesimal rigidity of the moduli
  spaces will also be briefly discussed. Similar to our previous survey
  articles \cite{lsy6, lsy5}, we will briefly describe the basic ideas of our proofs, the details of the proofs
  will be published soon, see \cite{lsy3, lsy4}.

 After introducing the definition of Weil-Petersson metric in
 Section \ref{wpme}, we discuss the fundamental curvature formula of Wolpert
 for the Weil-Petersson metric. For the reader's convenience we also briefly give a proof of the
 negativity of the Riemannian curvature of the Weil-Petersson
 metric. In Section \ref{rcprc} we discuss the Ricci and the perturbed Ricci
 metrics and their curvature formulas. In Section \ref{asymp} we describe the
 asymptotics of these metrics and their curvatures which are
 important for our understanding of their bounded geometry. In
 Section \ref{equivalence} we briefly discuss the equivalence of all of the
 complete metrics on Teichm\"uller spaces to the Ricci and the
 perturbed Ricci metrics, which is a simple corollary of our understanding
 of these two new metrics. In Section \ref{goodness} we discuss the goodness of the
 Weil-Petersson metric, the Ricci, the perturbed Ricci metric and
 the K\"ahler-Einstein metric. To prove the goodness we need much
 more
 subtle estimates on the connection and the curvatures of these
 metrics. Section \ref{natural} contains discussions of the dual Nakano
 negativity of the logarithmic tangent bundle of the moduli space
 and the naturalness of the Ricci and the perturbed Ricci metrics. In Section \ref{krf} we discuss the K\"ahler-Ricci flow and the \ke metric on the moduli space. There
 are many interesting corollaries from our understanding of the
 geometry of the moduli spaces. In Section \ref{app} we discuss the
 stability of the logarithmic cotangent bundle, the $L^2$
 cohomology and the infinitesimal rigidity of the moduli spaces as well as the Gauss-Bonnet theorem on the moduli space.

The \tei and moduli spaces of polarized Calabi-Yau (CY) manifolds
and Hyper-\ka manifolds are also important in mathematics and high
energy physics. In Section \ref{torelli} we will describe our recent
joint work with A. Todorov on the proof of global Torelli theorem of
the \tei space of polarized CY manifolds and Hyper-\ka manifolds. As
applications we will describe the construction of a global
holomorphic flat connection on the \tei space of CY manifolds and
the existence of \ke metrics on the Hodge completion of such \tei
spaces.

\section{The Weil-Petersson Metric and Its Curvature}\label{wpme}

Let $\M_g$ be the moduli space of Riemann surfaces of genus $g$
where $g\geq 2$. It is well known that the $\M_g$ is a complex
orbifold. The \tei space $\T_g$, as the space parameterizing marked
Riemann surfaces, is a smooth contractible pseudo-convex domain and
can be embedded into the Euclidean space of the same dimension.

\begin{remark}\label{obd}
Since $\M_g$ is only an orbifold, in the following when we work near
a point $p\in\M_g$ which is an orbifold point, we always work on a
local manifold cover of $\M_g$ around $p$. An alternative way is to
add a level structure on the moduli space so that it becomes smooth
\cite{shue1}. All the following results are still valid. In
particular, when we use the universal curve over the moduli space,
we always mean the universal curve over the local manifold cover.
When we deal with global properties of the moduli space, we can use
the moduli space with a level structure such that it becomes smooth.
We take quotients after we derive the estimates. We can also work on
the \tei space which is smooth.
\end{remark}

For any point $p\in\M_g$ we let $X_p$ be the corresponding Riemann
surface. By the Kodaira-Spencer theory we have the identification
\[
T_p^{1,0}\M_g\cong H^1\lb X_p,T_{X_p}^{1,0}\rb.
\]
It follows from Serre duality that
\[
\Omega_p^{1,0}\M_g\cong H^0\lb X_p,K_{X_p}^2\rb.
\]
By the Riemann-Roch theorem we know that the dimension
$\dim_\C\M_g=n=3g-3$.

The Weil-Petersson (WP) metric is the first known \ka metric on
$\M_g$. Ahlfors showed that the WP metric is \ka and its holomorphic
sectional curvature is bounded above by a negative constant which
only depends on the genus $g$. Royden conjectured that the Ricci
curvature of the WP metric is also bounded above by a negative
constant. This conjecture was proved by Wolpert \cite{wol3}.

Now we briefly describe the WP metric and its curvature formula.
Please see the works \cite{wol4}, \cite{wol5} of Wolpert for
detailed description and various aspects of the WP metric.

Let $\pi:\mathfrak X\to\M_g$ be the universal family over the moduli
space. For any point $s\in\M_g$ we let $X_s=\pi^{-1}(s)$ be the
corresponding smooth Riemann surface. Since the Euler characteristic
$\chi\lb X_s\rb=2-2g<0$, by the uniformization theorem we know that
each fiber $X_s$ is equipped with a unique \ke metric $\lambda$. In
the following we will always use the \ke metric $\lambda$ on $X_s$.
Let $z$ be any holomorphic coordinate on $X_s$. We have
\[
\pz\pzb\log\lambda=\lambda.
\]

Now we fix a point $s\in\M_g$ and let $(U,s_1,\cdots,s_n)$ be any
holomorphic coordinate chart on $\M_g$ around $s$. In the following
we will denote by $\pa_i$ and $\pz$ the local vector fields
$\frac{\pa}{\pa s_i}$ and $\frac{\pa}{\pa z}$ respectively. By the
Kodaira-Spencer theory and the Hodge theory we have the
identification
\[
T_s^{1,0}\M_g\cong \check H^1\lb X_s,T_{X_s}^{1,0}\rb \cong
\h^{0,1}\lb X_s,T_{X_s}^{1,0}\rb
\]
where the right side of the above formula is the space of harmonic
Beltrami differentials. In fact we can explicitly construct the
above identification. We let \[
a_i=a_i(z,s)=-\lambda^{-1}\pa_i\pzb\log\lambda
\]

and let
\[
v_i=\frac{\pa}{\pa s_i}+a_i\frac{\pa}{\pa z}.
\]
The vector field $v_i$ is a smooth vector field on $\pi^{-1}(U)$ and
is called the harmonic lift of $\frac{\pa}{\pa s_i}$. If we let
$B_i=\bar\pa_F v_i\in A^{0,1}\lb X_s,T_{X_s}^{1,0}\rb$ then $B_i$ is
harmonic and the map $\frac{\pa}{\pa s_i} \mapsto B_i$ is precisely
the Kodaira-Spencer map. Here $\pa_F$ is the operator in the fiber
direction. In local coordinates if we let $B_i=A_i d\bar z\otimes
\pz$ then $A_i=\pzb a_i$. Furthermore, it was proved by Schumacher
that if $\eta$ is any relative $(1,1)$-form on $\mathfrak X$ then
\begin{eqnarray}\label{20}
\frac{\pa}{\pa s_i}\int_{X_s}\eta=\int_{X_s} L_{v_i}\eta.
\end{eqnarray}

We note that although $A_i$ is a local smooth function on $X_s$, the
product
\[
A_i\bar A_j=B_i\cdot\bar B_j\in C^\infty(X_s)
\]
is globally defined. We let
\[
f_{i\bar j}=A_i\bar A_j\in C^\infty(X_s).
\]
The Weil-Petersson metric on $\M_g$ is given by
\[
h_{i\bar j}(s)=\int_{X_s}B_i\cdot\bar B_j\ dv= \int_{X_s}f_{i\bar
j}\ dv
\]
where $dv=\frac{\sqrt{-1}}{2}\lambda dz\wedge d\bar z$ is the volume
form on $X_s$ with respect to the \ke metric.

Now we describe the curvature formula of the WP metric. We let
$\Box=-\lambda^{-1}\pz\pzb$ be the Hodge-Laplacian acting on
$C^\infty(X_s)$. It is clear that the operator $\Box+1$ has no
kernel and thus is invertible. We let
\[
e_{i\bar j}=(\Box+1)^{-1}\lb f_{i\bar j}\rb\in C^\infty(X_s).
\]
The following curvature formula is due to Wolpert. See \cite{lsy1}
for the detailed proof.
\begin{proposition}\label{wpcurv}
Let $R_{i\bar jk\bar l}$ be the curvature of the WP metric. Then
\begin{eqnarray}\label{30}
R_{i\bar jk\bar l}=-\int_{X_s} \lb e_{i\bar j} f_{k\bar l}+ e_{i\bar
l} f_{k\bar j}\rb \ dv.
\end{eqnarray}
\end{proposition}

The curvature of the WP metric has very strong negativity property.
In fact we shall see in Section \ref{natural} that the WP metric is
dual Nakano negative. We collect the negativity property of the WP
metric in the following proposition.

\begin{proposition}\label{negativewp}
The bisectional curvature of the WP metric on the moduli space
$\M_g$ is negative. The holomorphic sectional and Ricci curvatures
of the WP metric are bounded above by negative constants.
Furthermore, the Riemannian sectional curvature of the WP metric is
also negative.
\end{proposition}

{\bf Proof.} These results are well known, see \cite{wol3}. Here we
give a short proof of the negativity of the Riemannian sectional
curvature of the WP metric for the reader's  convenience. The proof
follows from expressing the Riemannian sectional curvature in term
of complex curvature tensors and using the curvature formula
\eqref{30}.

In general, let $(X^n,g, J)$ be a \ka manifold. For any point $p\in
X$ and two orthonormal real tangent vectors $u,v\in T_p^\R X$, we
let $X=\frac{1}{2}\lb u-iJu\rb$ and $Y=\frac{1}{2}\lb v-iJv\rb$ and
we know that $X,Y\in T_p^{1,0}X$. We can choose holomorphic local
coordinate $s=(s_1,\cdots,s_n)$ around $p$ such that
$X=\frac{\pa}{\pa s_1}$. If $v=span_\R \{ u, Ju\}$, since $v$ is
orthogonal to $u$ and its length is $1$, we know $v=\pm Ju$. In this
case we have
\[
R(u,v,u,v)=R(u,Ju,u,Ju)=4R_{1\bar 1 1\bar 1}.
\]
Thus the Riemannian sectional curvature and the holomorphic
sectional curvature have the same sign.

If $v$ is not contained in the real plane spanned by $u$ and $Ju$ we
can choose the coordinate $s$ such that $X=\frac{\pa}{\pa s_1}$ and
$Y=\frac{\pa}{\pa s_2}$. In this case a direct computation shows
that
\begin{eqnarray}\label{rie10}
R(u,v,u,v)=2\lb R_{1\bar 1 2\bar 2}- Re \lb R_{1\bar 2 1\bar
2}\rb\rb.
\end{eqnarray}

Now we fix a point $p\in\M_g$ and let $u,v\in T_p^\R\M_g$. Let $X,Y$
be the corresponding $(1,0)$-vectors. Since we know that the
holomorphic sectional curvature of the WP metric is strictly
negative, we assume $v\notin span_\R\{ u,Ju\}$ and thus we can
choose holomorphic local coordinates $s=(s_1,\cdots,s_n)$ around $p$
such that $X=\frac{\pa}{\pa s_1}(p)$ and $Y=\frac{\pa}{\pa s_2}(p)$.
By formulas \eqref{rie10} and \eqref{30} we have
\begin{align}\label{rie40}
\begin{split}
R(u,v,u,v)=& -2\lb \int_{X_p} \lb e_{1\bar 1}f_{2\bar 2}+e_{1\bar 2}f_{2\bar 1}-2 Re ( e_{1\bar 2}f_{1\bar 2}) \rb dv \rb\\
=& -2\lb \int_{X_p} \lb e_{1\bar 1}f_{2\bar 2}+e_{1\bar 2}f_{2\bar
1}-e_{1\bar 2}f_{1\bar 2}-e_{2\bar 1}f_{2\bar 1}\rb dv\rb.
\end{split}
\end{align}

To prove the proposition we only need to show that
\begin{eqnarray}\label{rie50}
\int_{X_p} e_{1\bar 2}f_{2\bar 1}\ dv \geq \int_{X_p} Re \lb
e_{1\bar 2}f_{1\bar 2} \rb dv
\end{eqnarray}
and
\begin{eqnarray}\label{rie60}
\int_{X_p} e_{1\bar 1}f_{2\bar 2}\ dv\geq \int_{X_p} e_{1\bar
2}f_{2\bar 1}\ dv
\end{eqnarray}
and both equalities cannot hold simultaneously.

To prove inequality \eqref{rie50} we let $\alpha=Re (e_{1\bar 2})$
and $\beta=Im (e_{1\bar 2})$. Then we know
\[
\int_{X_p} e_{1\bar 2}f_{2\bar 1}\ dv=\int_{X_p} \lb \alpha
(\Box+1)\alpha +\beta (\Box+1)\beta \rb dv
\]
and
\[
\int_{X_p} Re \lb e_{1\bar 2}f_{1\bar 2} \rb dv= \int_{X_p} \lb
\alpha (\Box+1)\alpha -\beta (\Box+1)\beta \rb dv.
\]
Thus formula \eqref{rie50} reduces to
\[
 \int_{X_p} \beta (\Box+1)\beta\ dv\geq 0.
\]
However, we know
\[
\int_{X_p} \beta (\Box+1)\beta\ dv=\int_{X_p} \lb \Vert
\nabla'\beta\Vert^2+\beta^2 \rb dv\geq 0
\]
and the equality holds if and only if $\beta=0$. If this is the case
then we know that $e_{1\bar 2}$ is a real value function and
$f_{1\bar 2}$ is real valued too. Since $f_{1\bar 1}=A_1\bar A_1$
and $f_{1\bar 2}=A_1\bar A_2$ and $f_{1\bar 1}$ is real-valued we
know that there is a function $f\in C^\infty (X_p\setminus S, \R)$
such that $A_2=f(z)A_1$ on $X_p\setminus S$. Here $S$ is the set of
zeros of $A_1$. Since both $A_1$ and $A_2$ are harmonic, we know
that $\bar\pa^* A_1=\bar\pa^* A_2=0$. These reduce to $\pz(\lambda
A_1)=\pz(\lambda A_2)=0$ locally. It follows that $\pz
f\mid_{X_p\setminus S}=0$. Since $f$ is real-valued we know that $f$
must be a constant. But $A_1$ and $A_2$ are linearly independent
which is a contradiction. So the strict inequality \eqref{rie50}
always holds.

Now we prove formula \eqref{rie60}. Let $G(z,w)$ be the Green's
function of the operator $\Box+1$ and let $T=(\Box+1)^{-1}$. By the
maximum principle we know that $T$ maps positive functions to
positive functions. This implies that the Green's function $G$ is
nonnegative. Since $G(z,w)=G(w,z)$ is symmetric we know that
\begin{align}\label{rie70}
\begin{split}
\int_{X_p} e_{1\bar 1}f_{2\bar 2}\ dv=& \int_{X_p\times X_p}G(z,w) f_{1\bar 1}(w) f_{2\bar 2}(z)\ dv(w)dv(z)\\
=& \frac{1}{2} \int_{X_p\times X_p}G(z,w) \lb f_{1\bar 1}(w)
f_{2\bar 2}(z)+ f_{1\bar 1}(z) f_{2\bar 2}(w) \rb dv(w)dv(z).
\end{split}
\end{align}
Similarly we have
\begin{align}\label{rie80}
\begin{split}
\int_{X_p} e_{1\bar 2}f_{2\bar 1}\ dv=& \int_{X_p\times X_p}G(z,w) f_{1\bar 2}(w) f_{2\bar 1}(z)\ dv(w)dv(z)\\
=& \frac{1}{2} \int_{X_p\times X_p}G(z,w) \lb f_{1\bar 2}(w)
f_{2\bar 1}(z)+ f_{1\bar 2}(z) f_{2\bar 1}(w) \rb dv(w)dv(z).
\end{split}
\end{align}
Formula \eqref{rie60} follows from the fact that
\begin{align*}
\begin{split}
& f_{1\bar 1}(w) f_{2\bar 2}(z)+ f_{1\bar 1}(z) f_{2\bar 2}(w)-f_{1\bar 2}(w) f_{2\bar 1}(z)- f_{1\bar 2}(z) f_{2\bar 1}(w) \\
= & \left | A_1(z)A_2(w)-A_1(w)A_2(z)\right |^2\geq 0.
\end{split}
\end{align*}

\qed

Although the WP metric has very strong negativity properties, as we
shall see in Section \ref{asymp}, the WP metric is not complete and
its curvatures have no lower bound and this is very restrictive.

\section{The Ricci and Perturbed Ricci Metrics}\label{rcprc}

In \cite{lsy1} and \cite{lsy2} we studied two new \ka metrics: the
Ricci metric $\omega_\tau$ and the perturbed Ricci metric
$\omega_{\tilde\tau}$ on the moduli space $\M_g$. These new \ka
metrics are complete and have bounded geometry and thus have many
important applications. We now describe these new metrics.

Since the Ricci curvature of the WP metric has negative upper bound,
we define the Ricci metric
\[
\omega_\tau=-Ric\lb\omega_{_{WP}}\rb.
\]
We also define the perturbed Ricci metric to be a linear combination
of the Ricci metric and the WP metric
\[
\omega_{\tilde\tau}=\omega_\tau+C\omega_{_{WP}}
\]
where $C$ is a positive constants. In local coordinates we have
$\tau_{i\bar j}=-h^{k\bar l}R_{i\bar jk\bar l}$ and
$\tilde\tau_{i\bar j}=\tau_{i\bar j} +Ch_{i\bar j}$ where $R_{i\bar
jk\bar l}$ is the curvature of the WP metric.

Similar to curvature formula \eqref{30} of the WP metric we can
establish integral formulae for the curvature of the Ricci and
perturbed Ricci metrics. These curvature formulae are crucial in
estimating the asymptotics of these metrics and their curvature. To
establish these formulae, we need to introduce some operators. We
let
\[
P: C^\infty(X_s)\to A^{1,0}\lb T_{X_s}^{0,1}\rb
\]
be the operator defined by
\[
f\mapsto \pa \lb  \omega_{_{KE}}^{-1}\lrcorner\pa f\rb.
\]
In local coordinate we have $P(f)=\pz\lb \lambda^{-1}\pz f\rb
dz\otimes \pzb$. For each $1\leq k\leq n$ we let
\[
\xi_k:C^\infty(X_s)\to C^\infty(X_s)
\]
be the operator defined by
\[
f\mapsto \bar\pa^*\lb B_k\lrcorner \pa f\rb=-B_k \cdot P(f).
\]
In the local coordinate we have $\xi_k(f)=-\lambda^{-1}\pz\lb A_k
\pz f\rb$. Finally for any $1\leq k,l\leq n$ we define the operator
\[
Q_{k\bar l}:C^\infty(X_s)\to C^\infty(X_s)
\]
by
\[
Q_{k\bar l}(f)=\bar P\lb e_{k\bar l}\rb P(f)-2f_{k\bar l}\Box
f+\lambda^{-1}\pz f_{k\bar l}\pzb f.
\]
These operators are commutators of various classical operators on
$X_s$. See \cite{lsy1} for details. Now we recall the curvature
formulae of the Ricci and perturbed Ricci metrics established in
\cite{lsy1}. For convenience, we introduce the symmetrization
operator.
\begin{definition}
Let $U$ be any quantity which depends on indices $i,k,\alpha,\bar
j,\bar l, \bar\beta$. The symmetrization operator $\sigma_1$ is
defined by taking the summation of all orders of the triple
$(i,k,\alpha)$. That is
\begin{align*}
\begin{split}
\sigma_1(U(i,k,\alpha,\bar j,\bar l, \bar\beta))=& U(i,k,\alpha,\bar
j,\bar l, \bar\beta)+ U(i,\alpha,k,\bar j,\bar l, \bar\beta)+
U(k,i,\alpha,\bar j,\bar l, \bar\beta)\\ + & U(k,\alpha,i,\bar
j,\bar l, \bar\beta)+U(\alpha,i,k,\bar j,\bar l, \bar\beta)+
U(\alpha,k,i,\bar j,\bar l, \bar\beta).
\end{split}
\end{align*}
Similarly, $\sigma_2$ is the symmetrization operator of $\bar j$ and
$\bar \beta$ and $\widetilde{\sigma_1}$ is the symmetrization
operator of $\bar j$, $\bar l$ and $\bar \beta$.
\end{definition}

Now we can state the curvature formulae. We let $T=(\Box+1)^{-1}$ be
the operator in the fiber direction.
\begin{theorem}\label{riccicurv}
Let $s_1,\cdots,s_n$ be local holomorphic coordinates at $s \in
\M_g$ and let $\widetilde{R}_{i\bar j k\bar l}$ be the curvature of
the Ricci metric. Then at $s$, we have
\begin{align}\label{finalcurv}
\begin{split}
\widetilde{R}_{i\bar j k\bar l}=& -h^{\alpha\bar\beta}
\left\{\sigma_1\sigma_2\int_{X_s} \left\{T(\xi_k(e_{i\bar j}))
\bar{\xi}_l(e_{\alpha\bar\beta})+ T(\xi_k(e_{i\bar j}))
\bar{\xi}_\beta(e_{\alpha\bar l})
\right\}\ dv\right\}\\
&-h^{\alpha\bar\beta} \left\{\sigma_1\int_{X_s}Q_{k\bar l}(e_{i\bar
j}) e_{\alpha\bar\beta}\ dv
\right\}\\
&+\tau^{p\bar q}h^{\alpha\bar\beta}h^{\gamma\bar\delta}
\left\{\sigma_1\int_{X_s}\xi_k(e_{i\bar q}) e_{\alpha\bar\beta}\
dv\right\}\left\{ \widetilde\sigma_1\int_{X_s}\bar{\xi}_l(e_{p\bar
j})
e_{\gamma\bar\delta})\ dv\right\}\\
&+\tau_{p\bar j}h^{p\bar q}R_{i\bar q k\bar l}.
\end{split}
\end{align}
\end{theorem}

\begin{theorem}\label{perriccicurv}
Let $\tilde\tau_{i\bar j}=\tau_{i\bar j}+Ch_{i\bar j}$ where $\tau$
and $h$ are the Ricci and WP metrics respectively where $C>0$ is a
constant. Let $P_{i\bar j k\bar l}$ be the curvature of the
perturbed Ricci metric. Then we have
\begin{align}\label{finalpercurv}
\begin{split}
P_{i\bar j k\bar l}=&-h^{\alpha\bar\beta}
\left\{\sigma_1\sigma_2\int_{X_s} \left\{T(\xi_k(e_{i\bar j}))
\bar{\xi}_l(e_{\alpha\bar\beta})+ T (\xi_k(e_{i\bar j}))
\bar{\xi}_\beta(e_{\alpha\bar l})
\right\}\ dv\right\}\\
&-h^{\alpha\bar\beta} \left\{\sigma_1\int_{X_s}Q_{k\bar l}(e_{i\bar
j}) e_{\alpha\bar\beta}\ dv
\right\}\\
&+\widetilde\tau^{p\bar q}h^{\alpha\bar\beta}h^{\gamma\bar\delta}
\left\{\sigma_1\int_{X_s}\xi_k(e_{i\bar q}) e_{\alpha\bar\beta}\
dv\right\}\left\{ \widetilde\sigma_1\int_{X_s}\bar{\xi}_l(e_{p\bar
j})
e_{\gamma\bar\delta})\ dv\right\}\\
&+\tau_{p\bar j}h^{p\bar q}R_{i\bar q k\bar l} +CR_{i\bar j k\bar
l}.
\end{split}
\end{align}
\end{theorem}

In \cite{lsy1} and \cite{lsy2} we proved various properties of these
new metrics. Here we collect the important ones.
\begin{theorem}\label{rcprcpro}
The Ricci and perturbed Ricci metrics are complete \ka metrics on
$\M_g$. Furthermore we have
\begin{itemize}
\item These two metrics have bounded curvature.

\item The injectivity radius of the \tei space $\T_g$ equipped with any of these two metrics is bounded from below.

\item These metrics have \pg and thus the moduli space has finite volume when equipped with any of these metrics.

\item The perturbed Ricci metric has negatively pinched holomorphic sectional and Ricci curvatures when we choose the constant $C$ to be large enough.
\end{itemize}
\end{theorem}

The Ricci metric is also cohomologous to the \ke metric on $\M_g$ in
the sense of currents and hence can be used as the background metric
to estimate the \ke metric. We will discuss this in Section
\ref{equivalence}.

\section{Asymptotics}\label{asymp}

Since the moduli space $\M_g$ is noncompact, it is important to
understand the asymptotic behavior of the canonical metrics in order
to study their global properties. We first describe the local
pinching coordinates near the boundary of the moduli space by using
the plumbing construction of Wolpert.

Let $\mathcal M_g$ be the moduli space of Riemann surfaces of genus
$g \geq 2$ and let $\bar{\mathcal M}_g$ be its Deligne-Mumford
compactification \cite{dm1}. Each point $y \in \bar{\mathcal M}_g
\setminus \mathcal M_g$ corresponds to a stable nodal surface $X_y$.
A point $p \in X_y$ is a node if there is a neighborhood of $p$
which is isometric to the germ $\{ (u,v)\mid uv=0,\ |u|,|v|<1 \}
\subset \mathbb{C}^2$.

We first recall the rs-coordinate on a Riemann surface defined by
Wolpert in \cite{wol1}. There are two cases: the puncture case and
the short geodesic case. For the puncture case, we have a nodal
surface $X$ and a node $p\in X$. Let $a,b$ be two punctures which
are glued together to form $p$.
\begin{definition}
A local coordinate chart $(U,u)$ near $a$ is called rs-coordinate if
$u(a)=0$ where $u$ maps $U$ to the punctured disc $0<|u|<c$ with
$c>0$, and the restriction to $U$ of the K\"ahler-Einstein metric on
$X$ can be written as $\frac{1}{2|u|^2(\log |u|)^2} |du|^2$. The
rs-coordinate $(V,v)$ near $b$ is defined in a similar way.
\end{definition}
For the short geodesic case, we have a closed surface $X$, a closed
geodesic $\gamma \subset X$ with length $l <c_\ast$ where $c_\ast$
is the collar constant.
\begin{definition}
A local coordinate chart $(U,z)$ is called rs-coordinate at $\gamma$
if $\gamma \subset U$ where $z$ maps $U$ to the annulus
$c^{-1}|t|^{\frac{1}{2}}<|z| <c|t|^{\frac{1}{2}}$, and the
K\"ahler-Einstein metric on $X$ can be written as
$\frac{1}{2}(\frac{\pi}{\log |t|}\frac{1}{|z|}\csc \frac{\pi\log
|z|}{\log |t|})^2 |dz|^2$.
\end{definition}

By Keen's collar theorem \cite{ke1}, we have the following lemma:
\begin{lemma}\label{gcollar}
Let $X$ be a closed surface and let $\gamma$ be a closed geodesic on
$X$ such that the length $l$ of $\gamma$ satisfies $l <c_\ast$. Then
there is a collar $\Omega$ on $X$ with holomorphic coordinate $z$
defined on $\Omega$ such that
\begin{enumerate}
\item $z$ maps $\Omega$ to the annulus
$\frac{1}{c}e^{-\frac{2\pi^2}{l}}<|z|<c$ for $c>0$; \item the
K\"ahler-Einstein metric on $X$ restricted to $\Omega$ is given by
\begin{eqnarray}\label{precmetric}
(\frac{1}{2}u^2 r^{-2}\csc^2\tau) |dz|^2
\end{eqnarray}
where $u=\frac{l}{2\pi}$, $r=|z|$ and $\tau=u\log r$; \item the
geodesic $\gamma$ is given by the equation $|z|=
e^{-\frac{\pi^2}{l}}$.
\end{enumerate}
We call such a collar $\Omega$ a genuine collar.
\end{lemma}
We notice that the constant $c$ in the above lemma has a lower bound
such that the area of $\Omega$ is bounded from below. Also, the
coordinate $z$ in the above lemma is an rs-coordinate. In the
following, we will keep the notations $u$, $r$ and $\tau$.

Now we describe the local manifold cover of $\bar{\mathcal M}_g$
near the boundary. We take the construction of Wolpert \cite{wol1}.
Let $X_{0,0}$ be a stable nodal surface corresponding to a
codimension $m$ boundary point and let $p_1,\cdots,p_m$ be the nodes
of $X_{0.0}$. The smooth part $X_0=X_{0,0}\setminus \{
p_1,\cdots,p_m \}$ is a union of punctured Riemann surfaces. Fix the
rs-coordinate charts $(U_i,\eta_i)$ and $(V_i,\zeta_i)$ at $p_i$ for
$i=1,\cdots,m$ such that all the $U_i$ and $V_i$ are mutually
disjoint.

Now pick an open set $U_0 \subset X_0$ such that the intersection of
each connected component of $X_0$ and $U_0$ is a nonempty relatively
compact set and the intersection $U_0 \cap (U_i\cup V_i)$ is empty
for all $i$. We pick Beltrami differentials
$\nu_{m+1},\cdots,\nu_{n}$ which are supported in $U_0$ and span the
tangent space at $X_0$ of the deformation space of $X_0$. For
$s=(s_{m+1},\cdots,s_n)$, let $\nu(s)=\sum_{i=m+1}^n s_i\nu_i$. We
assume $|s|=(\sum |s_i|^2)^{\frac{1}{2}}$  small enough such that
$|\nu(s)|<1$. The nodal surface $X_{0,s}$ is obtained by solving the
Beltrami equation $\bar\partial w=\nu(s)\partial w$. Since $\nu(s)$
is supported in $U_0$, $(U_i,\eta_i)$ and $(V_i,\zeta_i)$ are still
holomorphic coordinates for $X_{0,s}$. Note that they are no longer
rs-coordinates. By the theory of Ahlfors and Bers \cite{ab1} and
Wolpert \cite{wol1} we can assume that there are constants
$\delta,c>0$ such that when $|s|<\delta$, $\eta_i$ and $\zeta_i$ are
holomorphic coordinates on $X_{0,s}$ with $0<|\eta_i|<c$ and
$0<|\zeta_i|<c$.

Now we assume $t=(t_1,\cdots,t_m)$ has small norm. We do the
plumbing construction on $X_{0,s}$ to obtain $X_{t,s}$ in the
following way. We remove from $X_{0,s}$ the discs $0<|\eta_i|\leq
\frac{|t_i|}{c}$ and $0<|\zeta_i|\leq \frac{|t_i|}{c}$ for each
$i=1,\cdots,m$, and identify $\frac{|t_i|}{c}<|\eta_i|< c$ with
$\frac{|t_i|}{c}<|\zeta_i|< c$ by the rule $\eta_i \zeta_i=t_i$.
This defines the surface $X_{t,s}$. The tuple
$(t_1,\cdots,t_m,s_{m+1},\cdots,s_n)$ are the local pinching
coordinates for the manifold cover of $\bar{\mathcal M}_g$. We call
the coordinates $\eta_i$ (or $\zeta_i$) the plumbing coordinates on
$X_{t,s}$ and the collar defined by $\frac{|t_i|}{c}<|\eta_i|< c$
the plumbing collar.
\begin{remark}
>From the estimate of Wolpert \cite{wol2}, \cite{wol1} on the length
of short geodesic, we have $u_i=\frac{l_i}{2\pi}\sim
-\frac{\pi}{\log|t_i|}$.
\end{remark}

Let $(t,s)=(t_1,\cdots,t_m,s_{m+1},\cdots,s_n)$ be the pinching
coordinates near $X_{0,0}$. For $|(t,s)|<\delta$, let $\Omega^j_c$
be the $j$-th genuine collar on $X_{t,s}$ which contains a short
geodesic $\gamma_j$ with length $l_j$. Let $u_j=\frac{l_j}{2\pi}$,
$u_0=\sum_{j=1}^m u_j+\sum_{j=m+1}^n |s_j|$, $r_j=|z_j|$ and
$\tau_j=u_j\log r_j$ where $z_j$ is the properly normalized
rs-coordinate on $\Omega^j_c$ such that
\[
\Omega^j_c=\{ z_j\mid c^{-1}e^{-\frac{2\pi^2}{l_j}}<|z_j|<c \}.
\]
>From the above argument, we know that the K\"ahler-Einstein metric
$\lambda$ on $X_{t,s}$, restrict to the collar $\Omega^j_c$, is
given by
\begin{eqnarray}\label{jmetric}
\lambda=\frac{1}{2}u_j^2 r_j^{-2}\csc^2\tau_j .
\end{eqnarray}
For convenience, we let $\Omega_c=\cup_{j=1}^m \Omega^j_c$ and
$R_c=X_{t,s}\setminus \Omega_c$. In the following, we may change the
constant $c$ finitely many times, clearly this will not affect the
estimates.

To estimate the WP, Ricci and perturbed Ricci metrics and their
curvatures, we first need to to find all the harmonic Beltrami
differentials $B_1,\cdots,B_n$ corresponding to the tangent vectors
$\frac{\partial}{\partial t_1},\cdots, \frac{\partial}{\partial
s_n}$. In \cite{ma1}, Masur constructed $3g-3$ regular holomorphic
quadratic differentials $\psi_1,\cdots,\psi_n$ on the plumbing
collars by using the plumbing coordinate $\eta_j$. These quadratic
differentials correspond to the cotangent vectors
$dt_1,\cdots,ds_n$.

However, it is more convenient to estimate the curvature if we use
the rs-coordinate on $X_{t,s}$ since we have the accurate form of
the K\"ahler-Einstein metric $\lambda$ in this coordinate. In
\cite{tr1}, Trapani used the graft metric constructed by Wolpert
\cite{wol1} to estimate the difference between the plumbing
coordinate and rs-coordinate and described the holomorphic quadratic
differentials constructed by Masur in the rs-coordinate. We collect
Trapani's results (Lemma 6.2-6.5, \cite{tr1}) in the following
theorem:
\begin{theorem}\label{imp}
Let $(t,s)$ be the pinching coordinates on $\bar{\mathcal M}_g$ near
$X_{0,0}$ which corresponds to a codimension $m$ boundary point of
$\bar{\mathcal M}_g$. Then there exist constants $M,\delta>0$ and
$1>c>0$ such that if $|(t,s)|<\delta$, then the $j$-th plumbing
collar on $X_{t,s}$ contains the genuine collar $\Omega^j_c$.
Furthermore, one can choose rs-coordinate $z_j$ on the collar
$\Omega_c^j$ such that the holomorphic quadratic differentials
$\psi_1,\cdots,\psi_n$ corresponding to the cotangent vectors
$dt_1,\cdots,ds_n$ have the form $\psi_i=\phi_i(z_j)dz_j^2$ on the
genuine collar $\Omega^j_c$ for $1\leq j \leq m$, where
\begin{enumerate}
\item $\phi_i(z_j)=\frac{1}{z_j^2}(q_i^j(z_j)+\beta_i^j)$ if
$i\geq m+1$; \item
$\phi_i(z_j)=(-\frac{t_j}{\pi})\frac{1}{z_j^2}(q_j(z_j)+\beta_j)$ if
$i=j$; \item $\phi_i(z_j)=(-\frac{t_i}{\pi})
\frac{1}{z_j^2}(q_i^j(z_j)+\beta_i^j)$ if $1\leq i \leq m$ and $i\ne
j$.
\end{enumerate}
Here $\beta_i^j$ and $\beta_j$ are functions of $(t,s)$, $q_i^j$ and
$q_j$ are functions of $(t,s,z_j)$ given by
\[
q_i^j(z_j)=\sum_{k<0}\alpha_{ik}^j(t,s)t_j^{-k}z_j^k
+\sum_{k>0}\alpha_{ik}^j(t,s)z_j^k
\]
and
\[
q_j(z_j)=\sum_{k<0}\alpha_{jk}(t,s)t_j^{-k}z_j^k
+\sum_{k>0}\alpha_{jk}(t,s)z_j^k
\]
such that
\begin{enumerate}
\item $\sum_{k<0}|\alpha_{ik}^j|c^{-k}\leq M$ and
$\sum_{k>0}|\alpha_{ik}^j|c^{k}\leq M$ if $i\ne j$; \item
$\sum_{k<0}|\alpha_{jk}|c^{-k}\leq M$ and
$\sum_{k>0}|\alpha_{jk}|c^{k}\leq M$; \item
$|\beta_i^j|=O(|t_j|^{\frac{1}{2}-\epsilon})$ with
$\epsilon<\frac{1}{2}$ if $i\ne j$; \item $|\beta_j|=(1+O(u_0))$.
\end{enumerate}
\end{theorem}

An immediate consequence is the precise asymptotics of the WP metric
which was computed in \cite{lsy1}. These asymptotic estimates were
also given by Wolpert in \cite{wol6}.
\begin{theorem}\label{wpasymp}
Let $(t,s)$ be the pinching coordinates and let $h$ be the WP
metric. Then
\begin{enumerate}
\item $h^{i\bar i}=2u_i^{-3}|t_i|^2(1+O(u_0))$ and $h_{i\bar i} =\frac{1}{2}\frac{u_i^{3}}{|t_i|^2}(1+O(u_0))$ for $1\leq i\leq m$;

\item $h^{i\bar j}=O(|t_it_j|)$ and $h_{i\bar j}=O\lb \frac{u_i^3u_j^3}{|t_it_j|}\rb$, if $1\leq i,j \leq m$ and $i\ne j$;

\item $h^{i\bar j}=O(1)$ and $h_{i\bar j}=O(1)$, if $m+1\leq i,j \leq n$;

\item $h^{i\bar j}=O(|t_i|)$ and $h_{i\bar j}=O\lb \frac{u_i^3}{|t_i|}\rb$ if $i\leq m < j$;

\item $h^{i\bar j}=O(|t_j|)$ and $h_{i\bar j}=O\lb \frac{u_j^3}{|t_j|}\rb$ if $j\leq m < i$.
\end{enumerate}
\end{theorem}

By using the asymptotics of the WP metric and the fact that
\[
B_i=\lambda^{-1}\sum_{j=1}^n h_{i\bar j}\bar\psi_j
\]
we can derive the expansion of the harmonic Beltrami differentials
corresponding to $\frac{\pa}{\pa t_i}$ and $\frac{\pa}{\pa s_j}$.

\begin{theorem}\label{aj10}
 For $c$ small, on the genuine collar $\Omega_c^j$, the coefficient functions
$A_i$ of the harmonic Beltrami differentials have the form:
\begin{enumerate}
\item $A_i=\frac{z_j}{\bar{z_j}}\sin^2\tau_j \lb
\bar{p_i^j(z_j)}+\bar{b_i^j}\rb$ if $i\ne j$;

\item
$A_j=\frac{z_j}{\bar{z_j}}\sin^2\tau_j\lb\bar{p_j(z_j)}+\bar{b_j}\rb$
\end{enumerate}
where
\begin{enumerate}
\item $p_i^j(z_j)=\sum_{k\leq -1}a_{ik}^j\rho_j^{-k}z_j^k
+\sum_{k\geq 1}a_{ik}^jz_j^k$ if $i\ne j$;

\item
$p_j(z_j)=\sum_{k\leq -1}a_{jk}\rho_j^{-k}z_j^k +\sum_{k\geq
1}a_{jk}z_j^k$.
\end{enumerate}
In the above expressions, $\rho_j=e^{-\frac{2\pi^2}{l_j}}$ and the
coefficients satisfy the following conditions:
\begin{enumerate}
\item $\sum_{k\leq -1}|a_{ik}^j|c^{-k}=O\lb u_j^{-2}\rb$ and
$\sum_{k\geq 1}|a_{ik}^j|c^{k}=O\lb u_j^{-2}\rb$\\ if $i\geq m+1$;

\item
$\sum_{k\leq -1}|a_{ik}^j|c^{-k}=O\lb \frac{u_i^3
u_j^{-2}}{|t_i|}\rb$ and $\sum_{k\geq
1}|a_{ik}^j|c^{k}=O\lb\frac{u_i^3u_j^{-2}}{|t_i|}\rb$ \\
if $i\leq m$ and $i\ne j$;

\item $\sum_{k\leq
-1}|a_{jk}|c^{-k}=O\lb \frac{u_j}{|t_j|}\rb$ and $\sum_{k\geq
1}|a_{jk}|c^{k}=O\lb\frac{u_j}{|t_j|}\rb$;

\item $|b_i^j|=O(u_j)$ if $i\geq m+1$;

\item $|b_i^j|=O\lb u_j\rb O\lb \frac{u_i^3}{|t_i|}\rb$ if $i\leq m$ and $i\ne j$;

\item
$b_j=-\frac{u_j}{\pi\bar{t_j}}(1+O(u_0))$.
\end{enumerate}
\end{theorem}

By a detailed study of the curvature of the WP metric we derived the
precise asymptotics of the Ricci metric in \cite{lsy1}.

\begin{theorem}\label{ricciest}
Let $(t,s)$ be the pinching coordinates. Then we have
\begin{enumerate}
\item $\tau_{i\bar i}=\frac{3}{4\pi^2}\frac{u_i^2}{|t_i|^2}(1+O(u_0))$ and $\tau^{i\bar i}=\frac{4\pi^2}{3}\frac{|t_i|^2}{u_i^2} (1+O(u_0))$, if $i \leq m$;

\item $\tau_{i\bar j}=O\bigg (\frac{u_i^2u_j^2}{|t_it_j|}(u_i+u_j)\bigg )$ and $\tau^{i\bar j}=O(|t_it_j|)$, if $i,j \leq m$ and $i\ne j$;

\item $\tau_{i\bar j}=O\lb \frac{u_i^2}{|t_i|}\rb$ and $\tau^{i\bar j}=O(|t_i|)$, if $i\leq m$ and $j\geq m+1$;

\item $\tau_{i\bar j}=O\lb \frac{u_j^2}{|t_j|}\rb$ and $\tau^{i\bar j}=O(|t_j|)$, if $j\leq m$ and $i\geq m+1$;

\item $\tau_{i\bar j}=O(1)$, if $i,j \geq m+1$.
\end{enumerate}
\end{theorem}

In \cite{lsy1} we also derived the asymptotics of the curvature of
the Ricci metric.

\begin{theorem}\label{mainholo}
Let $X_0\in\bar{\mathcal{M}_g}\setminus\mathcal{M}_g$ be a
codimension $m$ point and let $(t_1,\cdots,t_m,s_{m+1},\cdots,s_n)$
be the pinching coordinates at $X_0$ where $t_1,\cdots,t_m$
correspond to the degeneration directions. Then the holomorphic
sectional curvature is negative in the degeneration directions and
is bounded in the non-degeneration directions. More precisely, there
exists $\delta>0$ such that, if $|(t,s)|<\delta$, then
\begin{eqnarray}\label{important100}
\widetilde R_{i\bar ii\bar i}=-
\frac{3u_i^4}{8\pi^4|t_i|^4}(1+O(u_0))
\end{eqnarray}
if $i\leq m$ and
\begin{eqnarray}\label{important200}
\left | \widetilde R_{i\bar ii\bar i}\right |=O(1)
\end{eqnarray}
if $i\geq m+1$. Here $\tilde R$ is the curvature of the Ricci
metric.

Furthermore, on $\mathcal M_g$, the holomorphic sectional curvature,
the bisectional curvature and the Ricci curvature of the Ricci
metric are bounded from above and below.
\end{theorem}

In \cite{lsy3} and \cite{lsy4} we derived more precise estimates of
the curvature of the Ricci and perturbed Ricci metrics which we will
discuss in Section \ref{goodness}.

\section{Canonical Metrics and Equivalence}\label{equivalence}

In addition to the WP, Ricci and perturbed Ricci metrics on the
moduli space, there are several other canonical metrics on $\M_g$.
These include the \tei metric, the Kobayashi metric, the
Carath\'eodory metric, the \ke metric, the induced Bergman metric,
the McMullen metric and the asymptotic Poincar\'e metric.

Firstly, on any complex manifold there are two famous Finsler
metrics: the Carath\'eodory and Kobayashi metrics. Now we describe
these metrics.

Let $X$ be a complex manifold and of dimension $n$. let $\Delta_R$
be the disk in $\mathbb C$ with radius $R$. Let $\Delta=\Delta_1$
and let $\rho$ be the Poincar\'e metric on $\Delta$. Let $p\in X$ be
a point and let $v\in T_p X$ be a holomorphic tangent vector. Let
$\text{Hol}(X,\Delta_R)$ and $\text{Hol}(\Delta_R,X)$ be the spaces
of holomorphic maps from $X$ to $\Delta_R$ and from $\Delta_R$ to
$X$ respectively. The Carath\'eodory norm of the vector $v$ is
defined to be
\[
\Vert v\Vert_C=\sup_{f\in\text{Hol}(X,\Delta)}\Vert f_\ast
v\Vert_{\Delta,\rho}
\]
and the Kobayashi norm of $v$ is defined to be
\[
\Vert v\Vert_K=\inf_{f\in\text{Hol}(\Delta_R,X),\ f(0)=p,\
f'(0)=v}\frac{2}{R}.
\]
It is well known that the Carath\'eodory metric is bounded from
above by the Kobayashi metric after proper normalization. The first
known metric on the \tei space $\T_g$ is the \tei metric which is
also an Finsler metric. Royden showed that, on $\T_g$, the \tei
metric coincides with the Kobayashi metric. Generalizations and
proofs of Royden's theorem can be found in \cite{masar}.

Now we look at the \ka metrics. Firstly, since the \tei space $\T_g$
is a pseudo-convex domain, by the work of Cheng and Yau \cite{cy1}
and the later work of Yau, there exist a unique complete \ke metric
on $\T_g$ whose Ricci curvature is $-1$.

There is also a canonical Bergman metric on $\T_g$ which we describe
now. In general, let $X$ be any complex manifold, let $K_X$ be the
canonical bundle of $X$ and let $W$ be the space of $L^2$
holomorphic sections of $K_X$ in the sense that if $\sigma\in W$,
then
\[
\Vert\sigma\Vert_{L^2}^2=\int_X
(\sqrt{-1})^{n^2}\sigma\wedge\bar\sigma<\infty.
\]
The inner product on $W$ is defined to be
\[
(\sigma,\rho)=\int_X (\sqrt{-1})^{n^2}\sigma\wedge\bar\rho
\]
for all $\sigma,\rho\in W$. Let $\sigma_1,\sigma_2,\cdots$ be an
orthonormal basis of $W$. The Bergman kernel form is the
non-negative $(n,n)$-form
\[
B_X=\sum_{j=1}^\infty(\sqrt{-1})^{n^2}\sigma_j\wedge\bar\sigma_j.
\]

With a choice of local coordinates $z_i,\cdots,z_n$, we have
\[
B_X=BE_X(z,\bar z)(\sqrt{-1})^{n^2}dz_1\wedge\cdots\wedge dz_n
\wedge d\bar z_1\wedge\cdots\wedge d\bar z_n
\]
where $BE_X(z,\bar z)$ is called the Bergman kernel function. If the
Bergman kernel $B_X$ is positive, one can define the Bergman metric
\[
B_{i\bar j}=\frac{\partial^2\log BE_X(z,\bar z)}{\partial z_i
\partial \bar z_j}.
\]
The Bergman metric is well-defined and is nondegenerate if the
elements in $W$ separate points and the first jet of $X$.

It is easy to see that both the \ke metric and the Bergman metric on
the \tei space $\T_g$ are invariant under the action of the mapping
class group and thus descend down to the moduli space.
\begin{remark}\label{bergman}
We note that the induced Bergman metric on $\M_g$ is different from
the Bergman metric on $\M_g$.
\end{remark}

In \cite{mc} McMullen introduced another \ka metric  $g_{1/l}$ on
$\mathcal M_g$ which is equivalent to the Teichm\"uller metric. Let
$Log:\ \mathbb R_{+}\to [0,\infty)$ be a smooth function such that
\begin{enumerate}
\item $Log(x)=\log x$ if $x \geq 2$; \item $Log(x)=0$ if $x \leq
1$.
\end{enumerate}
For suitable choices of small constants $\delta,\epsilon>0$, the
K\"ahler form of the McMullen metric $g_{1/l}$ is
\[
\omega_{1/l}=\omega_{WP}-i\delta\sum_{l_{\gamma}(X)<\epsilon}\partial\bar\partial
Log\frac{\epsilon}{l_\gamma}
\]
where the sum is taken over primitive short geodesics $\gamma$ on
$X$.

Finally, since $\M_g$ is quasi-projective, there exists a
non-canonical asymptotic Poincar\'e metric $\omega_{_{P}}$ on
$\M_g$. In general, Let $\bar M$ be a compact projective manifold of
dimension $m$. Let $Y\subset \bar M$ be a divisor of normal
crossings and let $M=\bar M\setminus Y$. Cover $\bar M$ by
coordinate charts $U_1,\cdots,U_p,\cdots,U_q$ such that $(\bar
U_{p+1}\cup\cdots\cup\bar U_q)\cap Y=\emptyset$. We also assume
that, for each $1\leq \alpha \leq p$, there is a constant $n_\alpha$
such that $U_\alpha\setminus
Y=(\Delta^\ast)^{n_\alpha}\times\Delta^{m-n_\alpha}$ and on
$U_\alpha$, $Y$ is given by $z_1^\alpha\cdots
z_{n_\alpha}^\alpha=0$. Here $\Delta$ is the disk of radius
$\frac{1}{2}$ and $\Delta^\ast$ is the punctured disk of radius
$\frac{1}{2}$. Let $\{\eta_i\}_{1\leq i\leq q}$ be the partition of
unity subordinate to the cover $\{U_i\}_{1\leq i\leq q}$. Let
$\omega$ be a K\"ahler metric on $\bar M$ and let $C$ be a positive
constant. Then for $C$ large, the K\"ahler form
\[
\omega_{_{P}}=C\omega+\sum_{i=1}^p\sqrt{-1}\partial\bar\partial
\bigg (\eta_i\log\log\frac{1}{\left | z_1^i\cdots z_{n_i}^i\right |}
\bigg )
\]
defines a complete metric on $M$ with finite volume since on each
$U_i$ with $1\leq i\leq p$, $\omega{_{p}}$ is bounded from above and
below by the local Poincar\'e metric on $U_i$. We call this metric
the asymptotic Poincar\'e metric.

In 2004 we proved in \cite{lsy1} that all complete metrics on the
moduli space are equivalent. The proof is based on asymptotic
analysis of these metrics and Yau's Schwarz Lemma. It is an easy
corollaries of our understanding of the Ricci and the perturbed
Ricci metrics. In July 2004 we learned from the announcement of
S.-K. Yeung in Hong Kong University
 where he announced he could prove a small and easy part of our results about the equivalences
 of some of these metrics by using a bounded pluri-subharmonic function. We received a hard copy of
Yeung's paper in November 2004 where he used a method similar to
ours in \cite{lsy1} to compare the Bergman, the Kobayashi and the
Carath\'eodory metric. It should be interesting to see how one can
use the bounded psh function to derive these equivalences.

We recall that two metrics on $\M_g$ are equivalent if one metric is
bounded from above and below by positive constant multiples of the
other metric.
\begin{theorem}\label{eqall}
On the moduli space $\M_g$ the Ricci metric, the perturbed Ricci
metric, the \ke metric, the induced Bergman metric, the McMullen
metric, the asymptotic Poincar\'e metric, the Carath\'eodory metric
and the Teichm\"uller-Kobayashi metric are equivalent.
\end{theorem}

The equivalence of several of these metrics hold in more general
setting. In 2004 we defined the holomorphic homogeneous regular
manifolds in \cite{lsy1} which generalized the idea of Morrey.
\begin{definition}\label{hhrm}
A complex manifold $X$ of dimension $n$ is called holomorphic
homogeneous regular if there are positive constants $r<R$ such that
for each point $p\in X$ there is a holomorphic map $f_p:X\to \C^n$
which satisfies
\begin{enumerate}
\item $f_p(p)=0$;

\item $f_p:X\to f_p(X)$ is a biholomorphism;

\item $B_r\subset f_p(X)\subset B_R$ where $B_r$ and $B_R$ are Euclidean balls with center $0$ in $\C^n$.
\end{enumerate}
\end{definition}

In 2009 Yeung \cite{yeung} used the above definition without
appropriate reference which he called domain with uniform squeezing
property. It follows from the restriction properties of canonical
metrics and Yau's Schwarz Lemma that
\begin{theorem}\label{eqhhrm}
Let $X$ be a holomorphic homogeneous regular manifold. Then the
Kobayashi metric, the Bergman metric and the Carath\'eodory metric
on $X$ are equivalent.
\end{theorem}

\begin{remark}
It follows from the Bers embedding theorem that the \tei space of
genus $g$ Riemann surfaces is a holomorphic homogeneous regular
manifold if we choose $r=2$ and $R=6$ in Definition \ref{hhrm}.
\end{remark}

\section{Goodness of Canonical Metrics}\label{goodness}

In his work \cite{mum1}, Mumford defined the goodness condition to
study the currents of Chern forms defined by a singular Hermitian
metric on a holomorphic bundle over a quasi-projective manifold
where he generalized the Hirzebruch's proportionality theorem to
noncompact case. The goodness condition is a growth condition of the
Hermitian metric near the compactification divisor of the base
manifold. The major property of a good metric is that the currents
of its Chern forms define the Chern classes of the bundle. Namely
the Chern-Weil theory works in this noncompact case.

Beyond the case of homogeneous bundles over symmetric spaces
discussed by Mumford in \cite{mum1}, several natural bundles over
moduli spaces of Riemann surfaces give beautiful and useful
examples. In \cite{wol1}, Wolpert showed that the metric induced by
the hyperbolic metric on the relative dualizing sheaf over the
universal curve of moduli space of hyperbolic Riemann surfaces is
good. Later it was shown by Trapani \cite{tr1} that the metric
induced by the WP metric on the determinant line bundle of the
logarithmic cotangent bundle of the Deligne-Mumford moduli space is
good. In both cases, the bundles involved are line bundles in which
cases it is easier to estimate the connection and curvature. Other
than these, very few examples of natural good metrics are known.

The goodness of the WP metric has been a long standing open problem.
In this section we describe our work in \cite{lsy3} which gives a
positive answer to this problem.

We first recall the definition of good metrics and their basic
properties described in \cite{mum1}. Let $\bar X$ be a projective
manifold of complex dimension $\dim_\C \bar X=n$. Let $D\subset \bar
X$ be a divisor of normal crossing and let $X=\bar X\setminus D$ be
a Zariski open manifold. We let $\Delta_r$ be the open disk in $\C$
with radius $r$, let $\Delta=\Delta_1$,
$\Delta_r^*=\Delta_r\setminus\{0\}$ and $\Delta ^*=\Delta
\setminus\{0\}$. For each point $p\in D$ we can find a coordinate
chart $(U,z_1,\cdots,z_n)$ around $p$ in $\bar X$ such that
$U\cong\Delta^n$ and $V=U\cap X\cong \lb
\Delta^*\rb^{m}\times\Delta^{n-m}$. We assume that $U\cap D$ is
defined by the equation $z_1\cdots z_k=0$. We let $U(r)\cong
\Delta_r^n$ for $0<r<1$ and let $V(r)=U(r)\cap X$.

On the chart $V$ of $X$ we can define a local \po metric:
\begin{eqnarray}\label{localpo}
\omega_{loc}=\frac{\sqrt{-1}}{2}\sum_{i=1}^k \frac{dz_i\wedge d\bar
z_i}{2|z_i|^2 \lb \log|z_i|\rb^2}+
\frac{\sqrt{-1}}{2}\sum_{i=k+1}^n dz_i\wedge d\bar z_i.
\end{eqnarray}

Now we cover $D\subset \bar X$ by such coordinate charts
$U_1,\cdots,U_q$ and let $V_i=U_i\cap X$. We choose coordinates
$z_1^i,\cdots,z_n^i$ such that $D\cap U_i$ is given by $z_1^i\cdots
z_{m_i}^i=0$.

A \ka metric $\omega_g$ on $X$ has \pg if for each $1\leq i\leq q$
there are constants $0\leq r_i\leq 1$ and $0\leq c_i<C_i$ such that
$\omega_g\mid_{V_i(r_i)}$ is equivalent to the local \po metric
$\omega_{loc}^i$:
\[
c_i\omega_{loc}^i\leq \omega_g\mid_{V_i(r_i)} \leq
C_i\omega_{loc}^i.
\]

In \cite{mum1} Mumford defined differential forms with \nspg:
\begin{definition}\label{formpg}
Let $\eta\in A^p(X)$ be a smooth $p$-form. Then $\eta$ has \pg if
for each $1\leq i\leq q$ there exists a constant $c_i>0$ such that
for each point $s\in V_i\lb \frac{1}{2}\rb$ and tangent vectors
$t_1,\cdots,t_p\in T_s X$ one has
\[
\left | \eta(t_1,\cdots,t_p)\right |^2\leq c_i
\prod_{j=1}^p\omega_{loc}^i (t_j,t_j).
\]
The $p$-form $\eta$ is good if and only if both $\eta$ and $d\eta$
have \nspg.
\end{definition}

\begin{remark}
It is easy to see that the above definition does not depend on the
choice of the cover $(U_1,\cdots,U_q)$ but it does depend on the
compactification $\bar X$ of $X$.
\end{remark}

The above definition is local. We now give a global formulation.

\begin{lemma}\label{pglobal}
Let $\omega_g$ be a \ka metric on $X$ with \nspg. Then a $p$-form
$\eta\in A^p(X)$ has \pg if and only if $\Vert\eta\Vert_g<\infty$
where $\Vert\eta\Vert_g$ is the $C^0$ norm of $\eta$ with respect to
the metric $g$. Furthermore, the fact that $\eta$ has \pg is
independent of the choice of $g$. It follows that if $\eta_1\in
A^p(X)$ and $\eta_2\in A^q(X)$ have \nspg, then $\eta_1\wedge\eta_2$
also has \nspg.
\end{lemma}

Now we collect the basic properties of forms with \pg as described
in \cite{mum1}.
\begin{lemma}\label{pgpro}
Let $\eta\in A^p(X)$ be a form with \nspg. Then $\eta$ defines a
$p$-current on $\bar X$. Furthermore, if $\eta$ is good then
$d[\eta]=[d\eta]$.
\end{lemma}

Now we consider a holomorphic vector bundle $\bar E$ of rank $r$
over $\bar X$. Let $E=\bar E\mid_X$ and let $h$ be a Hermitian
metric on $E$. According to \cite{mum1} we have
\begin{definition}\label{metricgood}
The Hermitian metric $h$ is good if for any point $x\in D$, assume
$x\in U_i$ for some $i$, and any basis $e_1,\cdots,e_r$ of $\bar
E\mid_{U_i\lb \frac{1}{2}\rb}$, if we let
$h_{\alpha\bar\beta}=h(e_\alpha,e_\beta)$ then there exist positive
constants $c_i,d_i$ such that
\begin{enumerate}
\item $\left | h_{\alpha\bar\beta}\right |, \lb \det h\rb^{-1}\leq c_i \lb \sum_{j=1}^{m_i}\log|z_j|\rb^{2d_i}$;

\item the $1$-forms $\lb \pa h\cdot h^{-1}\rb_{\alpha\beta}$ are good on $V_i\lb \frac{1}{2}\rb$.
\end{enumerate}
\end{definition}

\begin{remark}\label{choose1}
A simple computation shows that the goodness of $h$ is independent
of the choice of the cover of $D$. Furthermore, to check whether a
metric $h$ is good or not  by using the above definition, we only
need to check the above two conditions for one choice of the basis
$e_1,\cdots,e_r$.
\end{remark}

The most important features of a good metric are
\begin{theorem}\label{megoodpro}
Let $h$ be a Hermitian metric on $E$. Then there is at most one
extension of $E$ to $\bar X$ for which $h$ is good. Furthermore, if
$h$ is a good metric on $E$, then the Chern forms $c_k(E,h)$ are
good and the current $[c_k(E,h)]=c_k(\bar E)\in H^{2k}(\bar X)$.
\end{theorem}
See \cite{mum1} for details. This theorem allows us to compute the
Chern classes by using Chern forms of a singular good metric.

Now we look at a special choice of the bundle $E$. In the following
we let $\bar E=T_{\bar X}(-\log D)$ to be the logarithmic tangent
bundle and let $E=\bar E\mid_X$. Let $U$ be one of the charts $U_i$
described above and assume $D\cap U$ is given by $z_1\cdots z_m=0$.
Let $V=V_i=U_i\cap X$. In this case a local frame of $\bar E$
restricting to $V$ is given by
\[
e_1=z_1\frac{\pa}{\pa z_1},\cdots,e_m=z_m\frac{\pa}{\pa z_m},\
e_{m+1}=\frac{\pa}{\pa z_{m+1}}, \cdots,e_n=\frac{\pa}{\pa z_n}.
\]
Let $g$ be any \ka metric on $X$. It induces a Hermitian metric
$\tilde g$ on $E$. In local coordinate $z=(z_1,\cdots,z_n)$ we have
\begin{eqnarray}\label{induce10}
\tilde g_{i\bar j}=
\begin{cases}
z_i\bar z_j g_{i\bar j}&\ \ \text{if}\ \ \ i,j\leq m\\
z_i g_{i\bar j}&\ \ \text{if}\ \ \ i\leq m<j\\
\bar z_j g_{i\bar j}&\ \ \text{if}\ \ \ j\leq m<i\\
g_{i\bar j}&\ \ \text{if}\ \ \ i,j>m.
\end{cases}
\end{eqnarray}

In the following we denote by $\pa_i$ the partial derivative
$\frac{\pa}{\pa z_i}$. Let
\[
\Gamma_{ik}^p=g^{p\bar q}\pa_i g_{k\bar q}
\]
be the Christoffel symbol of the \ka metric $g$ and let
\[
R_{ik\bar l}^p=g^{p\bar j}R_{i\bar jk\bar l}=g^{p\bar j}\lb
-\pa_k\pa_{\bar l}g_{i\bar j}+g^{s\bar t}\pa_k g_{i\bar t}\pa_{\bar
l}g_{s\bar j}\rb
\]
be the curvature of $g$. We define
\begin{eqnarray}\label{dik}
D_i^k=
\begin{cases}
\frac{z_i}{z_k} & \text{if}\ \ \ i,k\leq m\\
\frac{1}{z_k} & \text{if}\ \ \ k\leq m<i\\
z_i  & \text{if}\ \ \ i\leq m<k\\
1  & \text{if}\ \ \ i,k>m
\end{cases}
\end{eqnarray}
and we let
\begin{eqnarray}\label{gammaab}
\Lambda_i=
\begin{cases}
\frac{-1}{|z_i|\log |z_i|} & \text{if}\ \ \ i\leq m\\
1  & \text{if}\ \ \ i>m
\end{cases}.
\end{eqnarray}
Now we give an equivalent local condition of the metric $\tilde g$
on $E$ induced by the \ka metric $g$ to be good. We have
\begin{proposition}\label{logiff}
The metric $\tilde g$ on $E$ induced by $g$ is good on $V\lb
\frac{1}{2}\rb$ if and only if
\begin{align}\label{goodiff}
\begin{split}
& |\tilde g_{i\bar j}|,\ |z_1\cdots z_m|^{-2}\deg(g) \leq c\lb \sum_{i=1}^m \log|z_i|\rb^{2d} \ \ \text{for some constants}\ c,d>0\\
& \left | D_i^k \Gamma_{ip}^k \right |=O(\Lambda_p)\ \ \text{for all}\ 1\leq i,k,p\leq n\ \text{except}\ i=k=p\\
& \left | \frac{1}{t _i}+ \Gamma_{ii}^i \right |=O(\Lambda_i)\ \ \text{if}\ i\leq m\\
& \left | D_i^k R_{ip\bar q}^k \right |=O(\Lambda_p\Lambda_q).
\end{split}
\end{align}
\end{proposition}

In \cite{lsy3} we showed the goodness of the WP, Ricci and perturbed
Ricci metrics.
\begin{theorem}\label{goodmain}
Let $\M_g$ be the moduli space of genus $g$ Riemann surfaces. We
assume $g\geq 2$. Let $\bar\M_g$ be the Deligne-Mumford
compactification of $\M_g$ and let $D=\bar\M_g\setminus\M_g$ be the
compactification divisor which is a normal crossing divisor. Let
$\bar E=T_{\bar\M_g}(-\log D)$ and let $E=\bar E\mid_{\M_g}$. Let
$\hat h$, $\hat\tau$ and $\hat{\tilde\tau}$ be the metrics on $E$
induced by the WP, Ricci and perturbed Ricci metrics respectively.
Then $\hat h$, $\hat\tau$ and $\hat{\tilde\tau}$ are good in the
sense of Mumford.
\end{theorem}

This theorem is based on very accurate estimates of the connection
and curvature forms of these metrics. One of the difficulties is to
estimate the Gauss-Manin connection of the fiberwise \ke metric
where we use the compound graft metric construction of Wolpert
together with maximum principle.

\section{Negativity and Naturalness}\label{natural}

In Section \ref{wpme} we have seen various negative properties of
the WP metric. In fact, we showed in \cite{lsy3} that the WP metric
is dual Nakano negative. This means the complex curvature operator
of the dual metric of the WP metric is positive. We first recall the
precise definition of dual Nakano negativity of a Hermitian metric.

Let $(E,h)$ be a Hermitian holomorphic vector bundle of rank $m$
over a complex manifold $M$ of dimension $n$. Let $e_1,\cdots,e_m$
be a local holomorphic frame of $E$ and let $z_1,\cdots,z_n$ be
local holomorphic coordinates on $M$. The Hermitian metric $h$ has
expression $h_{i\bar j}=h\lb e_i,e_j\rb$ locally.

The curvature of $E$ is given by
\[
P_{i\bar j\alpha\bar\beta}=-\pa_\alpha\pa_{\bar\beta}h_{i\bar j}+
h^{p\bar q}\pa_\alpha h_{i\bar q}\pa_{\bar\beta}h_{p\bar j}.
\]

\begin{definition}\label{nadef}
The Hermitian vector bundle $(E,h)$ is Nakano positive if the
curvature $P$ defines a Hermitian metric on the bundle $E\otimes
T_M^{1,0}$. Namely, $P_{i\bar
j\alpha\bar\beta}C^{i\alpha}\bar{C^{j\beta}}>0$ for all $m\times n$
nonzero matrices $C$. The bundle $(E,h)$ is Nakano semi-positive if
$P_{i\bar j\alpha\bar\beta}C^{i\alpha}\bar{C^{j\beta}}\geq 0$. The
bundle is dual Nakano (semi-)negative if the dual bundle with dual
metric $(E^*,h^*)$ is Nakano (semi-)positive.
\end{definition}

We have proved the following theorem in \cite{lsy3}
\begin{theorem}\label{nakanomain}
Let $\M_g$ be the moduli space of Riemann surfaces of genus $g$
where $g\geq 2$. Let $h$ be the WP metric on $\M_g$. Then the
holomorphic tangent bundle $T^{1,0}\M_g$ equipped with the WP metric
$h$ is dual Nakano negative.
\end{theorem}

The dual Nakano negativity is the strongest negativity property of
the WP metric.

Now we look at the naturalness of the canonical metrics on the
moduli space. We let $\M_g$ be the moduli space of genus $g$ curves
where $g\geq 2$ and let $\bar\M_g$ be its Deligne-Mumford
compactification. We fix a point $p\in \bar\M_g\setminus\M_g$ of
codimension $m$ and let $X=X_p$ be the corresponding stable nodal
curve. The moduli space $\M(X)$ of the nodal surface $X$ is
naturally embedded into $\bar\M_g$. Furthermore, since each element
$Y$ in $\M(X)$ corresponds to a hyperbolic Riemann surface when we
remove the nodes from $Y$, the complement can be uniformized by the
upper half plane and thus there is a unique complete \ke metric on
$Y$ whose Ricci curvature is $-1$. We note that the moduli space
$\M(X)$ can be viewed as an irreducible component of the
intersection of $m$ compactification divisors.

By the discussion in Section \ref{wpme} there is a natural WP metric
$\hat h$ on $\M(X)$. The curvature formula \eqref{30} is still valid
for this WP metric and it is easy to see that the Ricci curvature of
the WP metric $\hat h$ is negative. We can take $\hat\tau=-Ric\lb
\omega_{\hat h}\rb$ to be the \ka form of a \ka metric on $\M(X)$.
This is the Ricci metric $\hat\tau$ on $\M(X)$.

In \cite{ma1} Masur showed that the WP metric $h$ on $\M_g$ extends
to $\bar\M_g$ and its restriction to $\M(X)$ via the natural
embedding $\M(X)\hookrightarrow\bar\M_g$ coincides with the WP
metric $\hat h$ on $\M(X)$. This implies the WP metric is natural.
In \cite{wol7} Wolpert showed that the WP Levi-Civita connection
restricted to directions which are almost tangential to the
compactification divisors limits to the lower dimensional WP
Levi-Civita connection. In \cite{lsy3} we proved the naturalness of
the Ricci metric.
\begin{theorem}\label{ricnatural}
The Ricci metric on $\M_g$ extends to $\bar\M_g$ in non-degenerating
directions. Furthermore, the restriction of the extension of $\tau$
to $\M(X)$ coincides with $\hat\tau$, the Ricci metric on $\M(X)$.
\end{theorem}

\section{The K\"ahler-Ricci Flow and K\"ahler-Einstein Metric on the Moduli Space}\label{krf}

The existence of the \ke metric on the \tei space was based on the
work of Cheng-Yau since the \tei space is pseudo-convex. By the
uniqueness we know that the \ke metric is invariant under the action
of the mapping class group and thus is also the \ke metric on the
moduli space. It follows from the later work of Yau that the \ke
metric is complete. However, the detailed properties of the \ke
metric remain unknown.

In \cite{lsy2} we proved the strongly bounded geometry property of
the \ke metric. We showed
\begin{theorem}\label{strong}
The \ke metric on the \tei space $\T_g$ has strongly bounded
geometry. Namely, the curvature and its covariant derivatives of the
\ke metric are bounded and the injectivity radius of the \ke metric
is bounded from below.
\end{theorem}

This theorem was proved in two steps. Firstly, we deform the Ricci
metric via the K\"ahler-Ricci flow
\begin{eqnarray}\label{bdd10}
\begin{cases}
\frac{\partial g_{i\bar j}}{\partial t}=-(R_{i\bar j}+g_{i\bar j})\\
g_{i\bar j}(0)=\tau_{i\bar j}
\end{cases}
\end{eqnarray}
Let $h=g(s)$ be the deformed metric at time $s\ll 1$. By the work of
Shi \cite{shi1} we know that the metric $h$ is equivalent to the
initial metric $\tau$ and is cohomologous to $\tau$ in the sense of
currents. Thus $h$ is complete and has \nspg. Furthermore, the
curvature and covariant derivatives of $h$ are bounded.

We then use the metric $h$ as a background metric to derive a priori
estimates for the \ke metric by using the Monge-Amper\'e equation
\[
\frac{\det\lb h_{i\bar j}+u_{i\bar j}\rb}{\det h_{i\bar j}}=e^{u+F}
\]
where $F$ is the Ricci potential of the metric $h$. If we denote by
$g$ the \ke metric and let
\[
S=g^{i\bar j}g^{k\bar l}g^{p\bar q}u_{;i\bar qk}u_{;\bar jp\bar l}
\]
and
\[
V=g^{i\bar j}g^{k\bar l}g^{p\bar q}g^{m\bar n} u_{;i\bar qk\bar
n}u_{;\bar jp\bar lm}+ g^{i\bar j}g^{k\bar l}g^{p\bar q}g^{m\bar n}
u_{;i\bar nkp}u_{;\bar jm\bar l\bar q}
\]
to be the third and fourth order quantities respectively. We have
\begin{align*}
\begin{split}
\Delta^{'}\left [(S+\kappa)V\right ]\geq & C_{1}\left
[(S+\kappa)V\right ]^2-C_{2} \left [(S+\kappa)V\right
]^{\frac{3}{2}}-C_{3}
\left [(S+\kappa)V\right ]\\
& -C_{4}\left [(S+\kappa)V\right ]^{\frac{1}{2}}
\end{split}
\end{align*}
where $\Delta'$ is the Laplace operator of the \ke metric $g$ and
$C_1>0$.

It follows from the mean value inequality that $S$ is bounded.
Furthermore, by the above estimate and the maximum principle we know
$V$ is bounded. In fact this method works for all higher order
derivatives of $u$ and we deduce that the \ke metric has strongly
bounded geometry.

The K\"ahler-Ricci flow and the goodness are closely tied together.
Firstly, since the most important feature of a Mumford good metric
is that the Chern-Weil theory still holds, we say metrics with this
property are intrinsic good. In \cite{lsy4} we showed
\begin{theorem}
Let $\bar X$ be a projective manifold with $\dim_\C \bar X=n$. Let
$D\subset\bar X$ be a divisor with normal crossings, let $X=\bar
X\setminus D$, let $\bar E=T_{\bar X}(-\log D)$ and let $E=\bar
E\mid_X$.

Let $\omega_g$ be a \ka metric on $X$ with bounded curvature and
\nspg. Assume $Ric(\omega_g)+\omega_g=\pa\bar\pa f$ where $f$ is a
bounded smooth function. Then
\begin{itemize}
\item There exists a unique \ke metric $\omega_{_{KE}}$ on $X$ with
\nspg.

\item The curvature and covariant derivatives of curvature of the
\ke metric are bounded.

\item If $\omega_g$ is intrinsic good, then $\omega_{_{KE}}$ is
intrinsic good. Furthermore, all metrics along the paths of
continuity and K\"ahler-Ricci flow are intrinsic good.
\end{itemize}

\end{theorem}

\section{Applications}\label{app}

In this last section we briefly look at some geometric applications
of the canonical metrics. The first application of the control of
the \ke metric is the stability of the logarithmic cotangent bundle
of the Deligne-Mumford moduli space. In \cite{lsy2} we proved
\begin{theorem}\label{stab}
Let $\bar E=T_{\bar\M_g}^*\lb \log D\rb$ be the logarithmic
cotangent bundle. Then $c_1(\bar E)$ is positive and $\bar E$ is
slope stable with respect to the polarization $c_1(\bar E)$.
\end{theorem}

An immediate consequence of the intrinsic goodness of the \ke metric
is the Chern number inequality. We have
\begin{theorem}
Let $\bar E= T_{\bar\M_g}(-\log D)$ be the logarithmic tangent
bundle of the moduli space. Then
\[
c_1(\bar E)^2 \leq \frac{6g-4}{3g-3}c_2(\bar E).
\]
\end{theorem}

An immediate consequence of the dual Nakano negativity and the
goodness of the WP metric is the positivity of the Chern numbers of
this bundle. We have
\begin{theorem}
The Chern numbers of the logarithmic cotangent bundle
$T_{\bar\M_g}^*(\log D)$ of the moduli spaces of Riemann surfaces
are all positive.
\end{theorem}

The dual Nakano negativity of a Hermitian metric on a bundle over a
compact manifold gives strong vanishing theorems by using Bochner
techniques. However, in our case the base variety $\M_g$ is only
quasi-projective. Thus we can only describe vanishing theorems of
the $L^2$ cohomology. In \cite{sap1}, Saper showed that the $L^2$
cohomology of the moduli space equipped with the Weil-Petersson
metric can be identified with the ordinary cohomology of the
Deligne-Mumford moduli space. Our situation is more subtle since the
natural object to be considered in our case is the tangent bundle
valued $L^2$ cohomology. Parallel to Saper's work, we proved in
\cite{lsy4}
\begin{theorem}\label{iden}
We have the following natural isomorphism
\[
H_{(2)}^* \lb \lb \M_g,\omega_\tau\rb,\lb
T_{\M_g},\omega_{_{WP}}\rb\rb \cong H^*\lb \bar\M_g, T_{\bar\M_g}\lb
-\log D\rb\rb.
\]
\end{theorem}

Now we combine the above result with the dual Nakano negativity of
the WP metric. In \cite{lsy4} we proved the following Nakano-type
vanishing theorem
\begin{theorem}\label{van}
The $L^2$ cohomology groups vanish:
\[
H_{(2)}^{0,q}\lb \lb \M_g,\omega_\tau\rb,\lb
T_{\M_g},\omega_{_{WP}}\rb\rb=0
\]
unless $q=3g-3$.
\end{theorem}

As a direct corollary we have
\begin{corollary}\label{rigid}
The pair $\lb \bar\M_g,D\rb$ is infinitesimally rigid.
\end{corollary}

Another important application of the properties of the Ricci,
perturbed Ricci and \ke metrics is the Gauss-Bonnet theorem on the
noncompact moduli space. Together with L. Ji, in \cite{jlsy} we
showed
\begin{theorem}
The Gauss-Bonnet theorem holds on the moduli space equipped with the
Ricci, perturbed Ricci or K\"ahler-Einstein metrics:
\[
\int_{\M_g} c_n(\omega_\tau) = \int_{\M_g} c_n(\omega_{\tilde\tau})
=\int_{\M_g} c_n(\omega_{_{KE}}) =\chi(\M_g) =
\frac{B_{2g}}{4g(g-1)}.
\]
Here $\chi(\M_g)$ is the orbifold Euler characteristic of $\M_g$ and
$n=3g-3$.
\end{theorem}

The explicit topological computation of the Euler characteristic of
the moduli space is due to Harer-Zagier \cite{hz1}. See also the
work of Penner \cite{penn1}.

As an application of the Mumford goodness of the WP metric and the
Ricci metric we have
\begin{theorem}
\[
\chi(T_{\bar\M_g}(-\log D))= \int_{\M_g} c_n(\omega_\tau)=
\int_{\M_g} c_n(\omega_{_{WP}})= \frac{B_{2g}}{4g(g-1)}
\]
where $n=3g-3$.
\end{theorem}

It is very hard to prove the Gauss-Bonnet theorem for the WP metric
directly since the WP metric is incomplete and its curvature is not
bounded. The proof is  based substantially on the Mumford goodness
of the WP metric.

By using the goodness of canonical metrics this theorem also gives
an explicit expression of the top log Chern number of the moduli
space.
\begin{theorem}
\[
\chi(\bar\M_g, T_{\bar\M_g}(-\log D))=\chi(\M_g)=
\frac{B_{2g}}{4g(g-1)}.
\]
\end{theorem}

\section{Global Torelli Theorem of the Teichm\"uller Spaces of Polarized Calabi-Yau Manifolds (Joint with Andrey Todorov)}\label{torelli}

The geometry of the \tei and moduli spaces of polarized Calabi-Yau
(CY) manifolds are the central objects in geometry and string
theory. One of the most important question in understanding the
geometry of the \tei
 and moduli space of polarized CY manifolds is the global Torelli problem which
asks whether the variation of
 polarized Hodge structures determines the marked polarized Calabi-Yau structure. In the rest of this article, after briefly discussing the deformation theory of CY manifolds and the geometry of period domain, we will describe the global Torelli theorem of the \tei spaces of polarized CY manifolds and its proof. See \cite{lsty1} for details.

Let $M$ be a Calabi-Yau manifold of dimension $\dim
_{\mathbb{C}}M=n$. Here we assume $n\geq 3$. Let $L$ be an ample
line bundle over $M$. By definition we assume that the canonical
bundle $K_M$ is trivial. Let $X$ be the underlying real
$2n$-dimensional manifold. We know that there is a nowhere vanishing
holomorphic $(n,0)$-form on $M$ which is unique up to scaling. The
\tei space of $(M,L)$ is the connected, simply connected, reduced
and irreducible manifold parameterizing triples $(M,L,(
\gamma_1,\cdots,\gamma_{b_n}))$ where $M$ is a CY manifold, $L$ is
the polarization and $(\gamma_1,\cdots,\gamma_{b_n})$ is a basis of
the middle homology group $H_n(X,{\mathbb{Z}})/tor$. Such triples
are called marked polarized CY manifolds.

\subsection{Deformation Theory of Polarized Calabi-Yau Manifolds}

We first recall the
 deformation of complex structures on a given smooth manifold. Let $X$ be a smooth
 manifold of dimension $\dim_{\mathbb{R}} X=2n$ and let $%
J_0$ be an integrable complex structure on $X$. We denote by
$M_0=(X,J_0)$ the corresponding complex manifold.

Let $\phi\in A^{0,1}\left (M_0,T_{M_0}^{1,0}\right )$ be a Beltrami
differential. We can view $\phi$ as a map
\begin{equation*}
\phi:\Omega^{1,0}(M_0)\to \Omega^{0,1}(M_0).
\end{equation*}
By using $\phi$ we define a new almost complex structure $J_\phi$ in
the following way. For a point $p\in M_0$ we pick a local
holomorphic coordinate chart $(U, z_1,\cdots,z_n)$ around $p$. Let
\begin{eqnarray}  \label{defmain}
\Omega_\phi^{1,0}(p)=\text{span}_{\mathbb{C}}\{ dz_1+\phi(dz_1),
\cdots, dz_n+\phi(dz_n)\}
\end{eqnarray}
and
\begin{equation*}
\Omega_\phi^{0,1}(p)=\text{span}_{\mathbb{C}}\{ d\bar
z_1+\bar\phi(d\bar z_1), \cdots, d\bar z_n+\bar\phi(d\bar z_n)\}
\end{equation*}
be the eigenspaces of $J_\phi$ with respect to the eigenvalue
$\sqrt{-1}$ and $-\sqrt{-1}$ respectively.

The almost complex structure $J_\phi$ is integrable if and only if
\begin{eqnarray}  \label{int}
\bar\partial\phi=\frac{1}{2}[\phi,\phi]
\end{eqnarray}
where $\bar\partial$ is the operator on $M_0$.

It was proved in \cite{tod1} and \cite{tian1} that the local
deformation of a polarized CY manifold is unobstructed.

\begin{theorem}
\label{unobs} The universal deformation space of a polarized CY
manifold is smooth.
\end{theorem}

The operation of contracting with $\Omega_0$ plays an important role
in converting bundle valued differential forms into ordinary
differential forms. The following lemma is the key step in the proof
of local Torelli theorem.

\begin{lemma}
\label{iso} Let $(M,L)$ be a polarized CY $n$-fold and let
$\omega_g$ be the unique CY metric in the class $[L]$. We pick a
nowhere vanishing holomorphic $(n,0)$-form $\Omega_0$ such that
\begin{eqnarray}  \label{normalization}
\left ( \frac{\sqrt{-1}}{2}\right )^n(-1)^{\frac{n(n-1)}{2}%
}\Omega_0\wedge\bar\Omega_0=\omega_g^n.
\end{eqnarray}
Then the map $\iota:A^{0,1}\left ( M, T_M^{1,0}\right )\to
A^{n-1,1}(M)$ given by $\iota(\phi)=\phi\lrcorner\Omega_0$ is an
isometry with respect to the natural Hermitian inner product on both
spaces induced by $\omega_g$. Furthermore, $\iota$ preserves the
Hodge decomposition.
\end{lemma}

In \cite{tod1} the existence of flat coordinates was established and
the flat coordinates played an important role in string theory
\cite{bcov1}. Here we recall this construction. Let $\mathfrak{X}$
be the universal family
over ${\mathcal{T}}$ and let $\pi$ be the projection map. For each $p\in{%
\mathcal{T}}$ we let $M_p=(X,J_p)$ be the corresponding CY manifold.
In the following we always use the unique CY metric on $M_p$ in the
polarization class $[L]$.

By the Kodaira-Spencer theory and Hodge theory, we have the
following identification
\begin{equation*}
T_p^{1,0}{\mathcal{T}}\cong {\mathbb{H}}^{0,1}\left (
M_p,T_{M_p}^{1,0}\right )
\end{equation*}
where we use ${\mathbb{H}}$ to denote the corresponding space of
harmonic forms. We have the following expansion of the Beltrami
differentials:

\begin{theorem}
\label{flatcoord} Let $\phi_1,\cdots,\phi_N \in
{\mathbb{H}}^{0,1}\left ( M_p,T_{M_p}^{1,0}\right )$ be a basis.
Then there is a unique power series
\begin{eqnarray}  \label{10}
\phi(\tau)=\sum_{i=1}^N \tau_i\phi_i +\sum_{|I|\geq 2}\tau^I\phi_I
\end{eqnarray}
which converges for $|\tau|<\varepsilon$. Here $I=(i_1,\cdots,i_N)$
is a multi-index, $\tau^I=\tau_1^{i_1}\cdots\tau_N^{i_N}$ and
$\phi_I\in A^{0,1}\left ( M_p,T_{M_p}^{1,0}\right )$. Furthermore,
if $\Omega$ is a nowhere vanishing holomorphic $(n,0)$-form, then
the family of Beltrami differentials $\phi(\tau)$ satisfy the
following conditions:
\begin{align}  \label{charflat}
\begin{split}
&\bar\partial_{M_p}\phi(\tau)=\frac{1}{2}[\phi(\tau),\phi(\tau)] \\
&\bar\partial_{M_p}^*\phi(\tau)=0 \\
&\phi_I\lrcorner\Omega =\partial_{M_p}\psi_I
\end{split}%
\end{align}
for each $|I|\geq 2$ where $\psi_I\in A^{n-2,1}(M_p)$. Furthermore,
by
shrinking $\varepsilon$ we can pick each $\psi_I$ appropriately such that $%
\sum_{|I|\geq 2}\tau^I\psi_I$ converges for $|\tau|<\varepsilon$.
\end{theorem}

The coordinates constructed in the above theorem are just the flat
coordinates described in \cite{bcov1}. They are unique up to affine
transformation and they are also the normal coordinates of the
Weil-Petersson metric at $p$. In fact the first equation in
\eqref{charflat} is the obstruction equation and the second equation
is the Kuranishi gauge which fixes the gauge in the fiber. The last
equation which characterized the flat coordinates around a point in
the \tei space is known as the Todorov gauge.

From Theorem \ref{flatcoord} the local Torelli theorem and the
Griffiths transversality follow immediately. However, Theorem
\ref{flatcoord} contains more information.

By using the local deformation theory, in \cite{tod1} Todorov
constructed a canonical local holomorphic section of the line bundle
$H^{n,0}=F^n$ over any flat coordinate chart $U\subset
{\mathcal{T}}$ in the form level. This canonical section plays a
crucial role in the proof of the global Torelli theorem.

We first consider the general construction of holomorphic
$(n,0)$-forms in \cite{tod1}.

\begin{lemma}
\label{constructn0} Let $M_0=(X,J_0)$ be a CY manifold where $J_0$
is the complex structure on $X$. Let $\phi\in A^{0,1}\left (
M_0,T_{M_0}^{1,0}\right )$ be a Beltrami differential on $M_0$ which
define an integrable complex structure $J_\phi$ and let
$M_\phi=(X,J_\phi)$ be the CY manifold whose underlying
differentiable manifold is $X$. Let $\Omega_0$ be a nowhere
vanishing holomorphic $(n,0)$-form on $M_0$ and let
\begin{eqnarray}  \label{n0expansion}
\Omega_\phi=\sum_{k=0}^n
\frac{1}{k!}(\wedge^k\phi\lrcorner\Omega_0).
\end{eqnarray}
Then $\Omega_\phi$ is a well-defined smooth $(n,0)$-form on
$M_\phi$. It is holomorphic with
respect to the complex structure $J_\phi$ if and only if $%
\partial(\phi\lrcorner\Omega_0)=0$. Here $\partial$ is the operator with
respect to the complex structure $J_0$.
\end{lemma}

By combining Lemma \ref{constructn0} and Theorem \ref{flatcoord} we
define the canonical family
\begin{eqnarray}  \label{can10}
\Omega^c=\Omega^c(\tau)=\sum_{k=0}^n \frac{1}{k!} \left (
\wedge^k\phi(\tau) \lrcorner\Omega_0\right )
\end{eqnarray}
and we have
\begin{corollary}
\label{expcoh} Let $\Omega^c(\tau)$ be a canonical family defined by %
\eqref{can10} where $\phi(\tau)$ is defined as in \eqref{charflat}.
Then we have the expansion
\begin{align}\label{cohexp10}
\begin{split}
[\Omega^c(\tau)]=& [\Omega_0]+\sum_{i=1}^N \tau_i[\phi_i\lrcorner\Omega_0]\\
&+\frac{1}{2}\sum_{i,j}\tau_i\tau_j \left [ {\mathbb{H}}(\phi_i\wedge\phi_j%
\lrcorner\Omega_0) \right ]+\Xi(\tau)
\end{split}
\end{align}
where $\Xi(\tau)\subset \bigoplus_{k=2}^n H^{n-k,k}(M_p)$ and $%
\Xi(\tau)=O(|\tau|^3)$.
\end{corollary}

The most important application of the cohomological expansion %
\eqref{cohexp10} is the invariance of the CY K\"ahler forms. The
theorem plays a central role in the proof of the global Torelli
theorem. This theorem was implicitly proved in \cite{bato1}. Please
see \cite{lsty1} for a simple and self-contained proof.

\begin{theorem}
\label{invarkaform} For each point $p\in{\mathcal{T}}$, let
$\omega_p$ be the K\"ahler form of the unique CY metric on $M_p$ in
the polarization class $[L]$. Then $\omega_p$ is invariant. Namely,
\begin{equation*}
\nabla^{GM}\omega_p=0.
\end{equation*}
Furthermore, since ${\mathcal{T}}$ is simply connected, we know that $%
\omega_p$ is a constant section of the trivial bundle
$A^2(X,{\mathbb{C}})$ over ${\mathcal{T}}$.
\end{theorem}

\subsection{The Teichm\"uller Space of Polarized Calabi-Yau Manifolds}

Now we recall the construction of the universal family of marked
polarized CY manifolds and the \tei space. See \cite{ltyz} for
details. Let $M$ be a CY manifold of dimension $\dim_{\mathbb{C}}
M=n\geq 3$. Let $L$ be an ample line bundle over $M$. We call a tuple $%
(M,L,\gamma_1,\cdots,\gamma_{h^n})$ a marked polarized CY manifold
if $M$ is a CY manifold, $L$ is a polarization of $M$ and $\left \{
\gamma_1,\cdots,\gamma_{h^n}\right \}$ is a basis of
$H_n(M,{\mathbb{Z}})/tor $.

\begin{remark}\label{simpli}
To simplify notations we assume in this section that a CY manifold
$M$ of dimension $n$ is simply connected and $h^{k,0}(M)=0$ for
$1\leq k\leq n-1$. All the results in this section hold when these
conditions are removed. This is due to the fact that we fix a
polarization.
\end{remark}

Since the Teichm\"uller space of $M$ with fixed marking and
polarization is constructed via GIT quotient, we need the following
results about group actions.

\begin{theorem}
\label{cons20} Let $(M,L,(\gamma_1,\cdots,\gamma_{b_n}))$ be a
marked polarized CY manifold and let $\pi:\mathfrak{X}\to
\mathcal{K}$ be the Kuranishi family of $M$. We let
$p\in\mathcal{K}$ such that $M=\pi^{-1}(p)$. If $G$ is a group of
holomorphic automorphisms of $M$ which preserve the
polarization $L$ and act trivially on $H_n(M,{\mathbb{Z}})$, then for any $%
q\in\mathcal{K}$ the group $G$ acts on $M_q=\pi^{-1}(q)$ as
holomorphic automorphisms.
\end{theorem}

Now we recall the construction of the Teichm\"uller space. We first
note that there is a constant $m_0>0$ which only depends on $n$ such
that for any polarized CY manifold $(M,L)$ of dimension $n$, the
line bundle $L^m$ is very ample for any $m\geq m_0$. We replace $L$
by $L^{m_0}$ and we still denote it by $L$.

Let $N_m=h^0(M,L^m)$. It follows from the Kodaira embedding theorem
that $M$ is embedded into $\mathbb{P}^{N_m-1}$ by the holomorphic
sections of $L^m$. Let $\mathcal{H}_L$ be the component of the
Hilbert scheme which contains $M$ and parameterizes smooth CY
varieties embedded in $\mathbb{P}^{N_m-1}$ with Hilbert polynomial
\begin{equation*}
P(m)=h^0\left ( M, L^m\right ).
\end{equation*}

We know that $\mathcal{H}_L$ is a smooth quasi-projective variety
and there
exists a universal family $\mathcal{X}_L\to \mathcal{H}_L$ of pairs $%
(M,(\sigma_0,\cdots,\sigma_{N_m}))$ where
$(\sigma_0,\cdots,\sigma_{N_m})$ is a basis of $H^0(M,L^m)$. Let
$\tilde{\mathcal{H}}_L$ be its universal
cover. By Theorem \ref{cons20} we know that the group $PGL\left ( N_m,{%
\mathbb{C}}\right )$ acts on $\tilde{\mathcal{H}}_L$ and the family $\tilde{%
\mathcal{X}}_L\to \tilde{\mathcal{H}}_L$ holomorphically and without
fixed points. Furthermore, it was proved in \cite{ltyz} and
\cite{ps} that the
group $PGL\left ( N_m,{\mathbb{C}}\right )$ also acts properly on $\tilde{%
\mathcal{H}}_L$.

We define the Teichm\"uller space of $M$ with polarization $L$ by
\begin{equation*}
{\mathcal{T}}={\mathcal{T}}_L(M)=\tilde{\mathcal{H}}_L/PGL\left ( N_m,{%
\mathbb{C}}\right ).
\end{equation*}

One of the most important features of the Teichm\"uller space is the
existence of universal family.

\begin{theorem}
\label{universal} There exist a family of marked polarized CY manifolds $\pi:%
\mathcal{U}_L\to{\mathcal{T}}_L(M)$ such that there is a point $p\in{%
\mathcal{T}}_L(M)$ with $M_p$ isomorphic to $M$ as marked polarized
CY manifolds and the family has the following properties:

\begin{enumerate}
\item ${\mathcal{T}}_L(M)$ is a smooth complex manifold of dimension $\dim_{%
\mathbb{C}}{\mathcal{T}}_L(M)=h^{n-1,1}(M)$.

\item For each point $q\in{\mathcal{T}}_L(M)$ there is a natural
identification
\begin{equation*}
T_q^{1,0}{\mathcal{T}}_L(M)\cong H^{0,1}\left (
M_q,T_{M_q}^{1,0}\right )
\end{equation*}
via the Kodaira-Spencer map.

\item Let $\rho:\mathcal{Y}\to \mathcal{C}$ be a family of marked polarized
CY manifold such that there is a point $x\in \mathcal{C}$ whose fiber $%
\rho^{-1}(x)$ is isomorphic to $M_p$ as marked polarized CY
manifolds. Then
there is a unique holomorphic map $f:(\mathcal{Y}\to \mathcal{C})\to (%
\mathcal{U}_L\to {\mathcal{T}}_L(M))$, defined up to biholomorphic
maps on the fibers whose induced maps on $H_n(M,{\mathbb{Z}})$ are
the identity map, such that $f$ maps the fiber $\rho^{-1}(x)$ to the
fiber $M_p$ and the
family $\mathcal{Y}$ is just the pullback of $\mathcal{U}_L$ via the map $f$%
. Furthermore, the map $\tilde f:\mathcal{C}\to{\mathcal{T}}_L(M)$
induced by $f$ is unique.
\end{enumerate}
\end{theorem}

It follows directly that

\begin{proposition}
\label{imp1} The Teichm\"uller space
${\mathcal{T}}={\mathcal{T}}_L(M)$ is a smooth complex manifold and
is simply connected.
\end{proposition}

\label{assumption} In the rest of this paper by the Teichm\"uller space ${%
\mathcal{T}}$ of $(M,L, (\gamma_1,\cdots,\gamma_{b_n}))$ we always
mean the
reduced irreducible component of $\tilde{\mathcal{H}}_L/PGL\left ( N_m,{%
\mathbb{C}}\right )$ with the fixed polarization $L$.

It follows from its construction and Theorem \ref{cons20} that the
universal family $\mathcal U_L$ over the \tei space $\T$ is
diffeomorphic to $M_p\times\T$ as a $C^\infty$ family where $p\in\T$
is any point and $M_p$ is the corresponding CY manifold.

\subsection{The Classifying Space of Variation of Polarized Hodge Structures}

Now we recall the construction of the classifying space of variation
of polarized Hodge structures and its basic properties such as the
description of its real and complex tangent spaces and the Hodge
metric. See \cite{schmid1} for details.

In the construction of the Teichm\"uller space ${\mathcal{T}}$ we
fixed a
marking of the background manifold $X$, namely a basis of $H_n(X,{\mathbb{Z}}%
)/tor$. This gives us canonical identifications of the middle
dimensional de Rahm cohomology of different fibers over
${\mathcal{T}}$. Namely for any two distinct points
$p,q\in{\mathcal{T}}$ we have the canonical identification
\begin{equation*}
H^n(M_p)\cong H^n(M_q)\cong H^n(X)
\end{equation*}
where the coefficient ring is ${\mathbb{Q}}$, ${\mathbb{R}}$ or
${\mathbb{C}} $.

Since the polarization $[L]$ is an integeral class, it defines a map
\begin{equation*}
L:H^n(X,{\mathbb{Q}})\to H^{n+2}(X,{\mathbb{Q}})
\end{equation*}
given by $A\mapsto c_1(L)\wedge A$ for any $A\in
H^n(X,{\mathbb{Q}})$. We denote by $H_{pr}^n(X)=\ker(L)$ the
primitive cohomology groups where,
again, the coefficient ring is ${\mathbb{Q}}$, ${\mathbb{R}}$ or ${\mathbb{C}%
}$. For any $p\in{\mathcal{T}}$ we let $H_{pr}^{k,n-k}(M_p)=H^{k,n-k}(M_p)%
\cap H_{pr}^n(M_p,{\mathbb{C}})$ and denote its dimension by
$h^{k,n-k}$. The Poincar\'e bilinear form $Q$ on
$H_{pr}^n(X,{\mathbb{Q}})$ is defined by
\begin{equation*}
Q(u,v)=(-1)^{\frac{n(n-1)}{2}}\int_X u\wedge v
\end{equation*}
for any $d$-closed $n$-forms $u,v$ on $X$. The bilinear form $Q$ is
symmetric if $n$ is even and is skew-symmetric if $n$ is odd.
Furthermore, $Q $ is non-degenerate and can be extended to
$H_{pr}^n(X,{\mathbb{C}})$ bilinearly. For any point
$q\in{\mathcal{T}}$ we have the Hodge decomposition
\begin{eqnarray}  \label{cl10}
H_{pr}^n(M_q,{\mathbb{C}})=H_{pr}^{n,0}(M_q,{\mathbb{C}})\oplus\cdots\oplus
H_{pr}^{0,n}(M_q,{\mathbb{C}})
\end{eqnarray}
which satisfies
\begin{eqnarray}  \label{cl20}
\dim_{\mathbb{C}} H_{pr}^{k,n-k}(M_q,{\mathbb{C}})=h^{k,n-k}
\end{eqnarray}
and the Hodge-Riemann relations
\begin{eqnarray}  \label{cl30}
Q\left ( H_{pr}^{k,n-k}(M_q,{\mathbb{C}}), H_{pr}^{l,n-l}(M_q,{\mathbb{C}}%
)\right )=0\ \ \text{unless}\ \ k+l=n
\end{eqnarray}
and
\begin{eqnarray}  \label{cl40}
\left (\sqrt{-1}\right )^{2k-n}Q\left ( v,\bar v\right )>0\ \
\text{for}\ \ v\in H_{pr}^{k,n-k}(M_q,{\mathbb{C}})\setminus\{0\}.
\end{eqnarray}

The above Hodge decomposition of $H_{pr}^n(M_q,{\mathbb{C}})$ can
also be described via the Hodge filtration. Let $f^k=\sum_{i=k}^n
h^{i,n-i}$. We let
\begin{equation*}
F^k=F^k(M_q)=H_{pr}^{n,0}(M_q,{\mathbb{C}})\oplus\cdots\oplus
H_{pr}^{k,n-k}(M_q,{\mathbb{C}})
\end{equation*}
and we have decreasing filtration
\begin{equation*}
H_{pr}^n(M_q,{\mathbb{C}})=F^0(M_q)\supset\cdots\supset F^n(M_q).
\end{equation*}
We know that
\begin{eqnarray}  \label{cl45}
\dim_{\mathbb{C}} F^k=f^k,
\end{eqnarray}
\begin{eqnarray}  \label{cl46}
H_{pr}(X,{\mathbb{C}})=F^{k}(q)\oplus \bar{F^{n-k+1}(q)}
\end{eqnarray}
and
\begin{eqnarray}  \label{cl48}
H_{pr}^{k,n-k}(M_q,{\mathbb{C}})=F^k(M_q)\cap\bar{F^{n-k}(M_q)}.
\end{eqnarray}
In term of the Hodge filtration $F^n\subset\cdots\subset F^0=H_{pr}^n(M_q,{%
\mathbb{C}})$ the Hodge-Riemann relations can be written as
\begin{eqnarray}  \label{cl50}
Q\left ( F^k,F^{n-k+1}\right )=0
\end{eqnarray}
and
\begin{eqnarray}  \label{cl60}
Q\left ( Cv,\bar v\right )>0 \ \ \text{if}\ \ v\ne 0
\end{eqnarray}
where $C$ is the Weil operator given by $Cv=\left (\sqrt{-1}\right
)^{2k-n}v$ when $v\in H_{pr}^{k,n-k}(M_q,{\mathbb{C}})$. The
classifying space $D$ of variation of polarized Hodge structures
with data \eqref{cl45} is the space of all such Hodge filtrations
\begin{equation*}
D=\left \{ F^n\subset\cdots\subset F^0=H_{pr}^n(X,{\mathbb{C}})\mid %
\eqref{cl45}, \eqref{cl50} \text{ and } \eqref{cl60} \text{ hold}
\right \}.
\end{equation*}
The compact dual $\check D$ of $D$ is
\begin{equation*}
\check D=\left \{ F^n\subset\cdots\subset F^0=H_{pr}^n(X,{\mathbb{C}})\mid %
\eqref{cl45} \text{ and } \eqref{cl50} \text{ hold} \right \}.
\end{equation*}
The classifying space $D\subset \check D$ is an open set. We note
that the
conditions \eqref{cl45}, \eqref{cl50} and \eqref{cl60} imply the identity %
\eqref{cl46}.

An important feature of the variation of polarized Hodge structures
is that both $D$ and $\check D$ can be written as quotients of
semi-simple Lie
groups. Let $H_{\mathbb{R}}=H_{pr}^n(X,{\mathbb{R}})$ and $H_{\mathbb{C}}%
=H_{pr}^n(X,{\mathbb{C}})$. We consider the real and complex
semi-simple Lie groups
\begin{equation*}
G_{\mathbb{R}}=\{\sigma\in GL(H_{\mathbb{R}})\mid Q(\sigma u, \sigma
v)=Q(u,v)\}
\end{equation*}
and
\begin{equation*}
G_{\mathbb{C}}=\{\sigma\in GL(H_{\mathbb{C}})\mid Q(\sigma u, \sigma
v)=Q(u,v)\}.
\end{equation*}
The real group $G_{\mathbb{R}}$ acts on $D$ and the complex group $G_{%
\mathbb{C}}$ acts on $\check D$ where both actions are transitive.
This implies that both $D$ and $\check D$ are smooth. Furthermore,
we can embed the real group into the complex group naturally as real
points.

We now fix a reference point $O=\{F_0^k\}\in D\subset\check D$ and
let $B$ be the isotropy group of $O$ under the action of
$G_{\mathbb{C}}$ on $\check D$. Let $\left\{ H_0^{k,n-k}\right\}$ be
the corresponding Hodge
decomposition where $H_0^{k,n-k}=F_0^k\cap \bar{F_0^{n-k}}$. Let $V=G_{%
\mathbb{R}}\cap B$. Then we have
\begin{eqnarray}  \label{cl70}
D=G_{\mathbb{R}}/V \ \ \ \text{and}\ \ \ \check D=G_{\mathbb{C}}/B.
\end{eqnarray}
Following the argument in \cite{schmid1} we let
\begin{equation*}
H_0^{+}=\bigoplus_{i\text{ is even}}H_0^{i,n-i} \ \ \ \ H_0^{-}=\bigoplus_{i%
\text{ is odd}}H_0^{i,n-i}
\end{equation*}
and let $K$ be the isotropy group of $H_0^{+}$ in $G_{\mathbb{R}}$.
We note that $H_0^+$ and $H_0^-$ are defined over ${\mathbb{R}}$ and
are orthogonal with respect to $Q$ when $n$ is even. When $n$ is odd
they are conjugate to each other. Thus $K$ is also the isotropy
group of $H_0^-$. In both cases $K$ is the maximal compact subgroup
of $G_{\mathbb{R}}$ containing $V$. This implies that $\tilde
D=G_{\mathbb{R}}/K$ is a symmetric space of noncompact
type and $D$ is a fibration over $\tilde D$ whose fibers are isomorphic to $%
K/V$.

\begin{remark}
\label{question} In the following we will only consider primitive
cohomology classes and we will drop the mark $``$pr$"$. Furthermore,
Since we only need to use the component of $G_{\mathbb{R}}$
containing the identity, we will denote again by $G_{\mathbb{R}}$
and $K$ the components of the real group and its corresponding
maximal compact subgroup which contain the identity.
\end{remark}

We fix a point $p\in{\mathcal{T}}$ and let $O=\Phi(p)\in
D\subset\check D$.
For $0\leq k\leq n$ we let $H_0^{k,n-k}=H^{k,n-k}(M_p)$. Now we let ${%
\mathfrak{g}}={\mathfrak{g}}_{\mathbb{C}}$ be the Lie algebra of $G_{\mathbb{%
C}}$ and let ${\mathfrak{g}}_0={\mathfrak{g}}_{\mathbb{R}}$ be the
Lie algebra of $G_{\mathbb{R}}$. The real Lie algebra
${\mathfrak{g}}_0$ can also be embedded into ${\mathfrak{g}}$
naturally as real points. The Hodge structure $\left\{
H_0^{k,n-k}\right \}$ induces a weight $0$ Hodge
structure on ${\mathfrak{g}}$. Namely ${\mathfrak{g}}=\bigoplus_p {\mathfrak{%
g}}^{p,-p}$ where
\begin{equation*}
{\mathfrak{g}}^{p,-p}=\left\{ X\in {\mathfrak{g}}\mid X\left
(H_0^{k,n-k}\right )\subset H_0^{k+p,n-k-p}\right \}.
\end{equation*}

Let $B$ be the isotropy group of $O\in\check D$ under the action of $G_{%
\mathbb{C}}$ and let $\mathfrak{b}$ be the Lie algebra of $B$. Then
\begin{equation*}
\mathfrak{b}=\bigoplus_{p\geq 0}{\mathfrak{g}}^{p,-p}.
\end{equation*}

Let $V=B\cap G_{\mathbb{R}}$ be the isotropy group of $O\in D$ under
the action of $G_{\mathbb{R}}$ and let $\mathfrak{v}$ be its Lie
algebra. We have
\begin{equation*}
\mathfrak{v}=\mathfrak{b}\cap {\mathfrak{g}}_0\subset
{\mathfrak{g}}.
\end{equation*}
Now we have
\begin{equation*}
\mathfrak{v}={\mathfrak{g}}_0\cap\mathfrak{b}={\mathfrak{g}}_0\cap\mathfrak{b%
}\cap\bar{\mathfrak{b}}={\mathfrak{g}}_0\cap {\mathfrak{g}}^{0,0}.
\end{equation*}

Let $\theta$ be the Weil operator of the weight $0$ Hodge structure on ${%
\mathfrak{g}}$. Then for any $v\in {\mathfrak{g}}^{p,-p}$ we have $%
\theta(v)=(-1)^pv$. The eigenvalues of $\theta$ are $\pm 1$. Let ${\mathfrak{%
g}}^+$ be the eigenspace of $1$ and let ${\mathfrak{g}}^-$ be the
eigenspace of $-1$. Then we have
\begin{equation*}
{\mathfrak{g}}^+=\bigoplus_{p \text{ even}}{\mathfrak{g}}^{p,-p} \ \
\ \text{ and }\ \ \ {\mathfrak{g}}^-=\bigoplus_{p \text{
odd}}{\mathfrak{g}}^{p,-p}.
\end{equation*}
We note here that, in the above expression, $p$ can be either
positive or negative. Let $\mathfrak{k}$ be the Lie algebra of $K$,
the maximal compact
subgroup of $G_{\mathbb{R}}$ containing $V$. By the work of Schmid \cite%
{schmid1} we know that

\begin{lemma}
\label{cartan} The Lie algebra $\mathfrak{k}$ is given by $\mathfrak{k}={%
\mathfrak{g}}_0\cap {\mathfrak{g}}^+$. Furthermore, if we let $\mathfrak{p}%
_0={\mathfrak{g}}_0\cap {\mathfrak{g}}^-$, then
\begin{equation*}
{\mathfrak{g}}_0=\mathfrak{k}\oplus \mathfrak{p}_0
\end{equation*}
is a Cartan decomposition of ${\mathfrak{g}}_0$. The space
$\mathfrak{p}_0$ is $Ad_V$ invariant.
\end{lemma}

We call such a Cartan decomposition the canonical Cartan
decomposition. Here we recall that if
${\mathfrak{g}}_0={\mathfrak{k}}\oplus {\mathfrak{p}}_0$
is a Cartan decomposition of the real semisimple Lie algebra ${\mathfrak{g}}%
_0$, then we know that ${\mathfrak{k}}$ is a Lie subalgebra, $[{\mathfrak{p}}%
_0,{\mathfrak{p}}_0]\subset {\mathfrak{k}}$ and $[{\mathfrak{p}}_0,{%
\mathfrak{k}}]\subset {\mathfrak{p}}_0$.

By the expression of $\mathfrak{v}$ and $\mathfrak{k}$ we have the
identification
\begin{eqnarray}  \label{kv100}
\mathfrak{k}/\mathfrak{v}\cong {\mathfrak{g}}_0\cap \left
(\bigoplus_{p\ne 0, \ p \text{ is even}} {\mathfrak{g}}^{p,-p}\right
)
\end{eqnarray}
and the identification
\begin{eqnarray}  \label{killing10}
T_O^{\mathbb{R}} D\cong \mathfrak{k}/\mathfrak{v}\oplus
\mathfrak{p}_0.
\end{eqnarray}

Now we look at the complex structures on $D$. By the above
identification we know that for each element $X\in T_O^{\mathbb{R}}
D$ we have the unique decomposition $X=X_++X_-$ where
$X_+\in\bigoplus_{p>0}{\mathfrak{g}}^{-p,p}$ and $X_-
\in\bigoplus_{p>0}{\mathfrak{g}}^{p,-p}$. We define the complex
structure $J$ on $T_O^{\mathbb{R}} D$ by
\begin{eqnarray}  \label{cxstructured}
JX=iX_+-iX_-.
\end{eqnarray}
Now we use left translation by elements in $G_{\mathbb{R}}$ to move
this complex structure to every point in $D$. Namely, for any point
$\alpha\in D$
we pick $g\in G_{\mathbb{R}}$ such that $g(O)=\alpha$. If $X\in T_\alpha^{%
\mathbb{R}} D$, then we define $JX=\left ( l_g\right )_*\circ J
\circ \left ( l_{g^{-1}}\right )_*(X)$.

\begin{lemma}
\label{cxstru} $J$ is an invariant integrable complex structure on
$D$. Furthermore, it coincides with the complex structure on $D$
induced by the inclusion $D\subset\check D=G_{\mathbb{C}}/B$.
\end{lemma}

This lemma is well known. See \cite{grsch1} and \cite{ls1} for
details.

There is a natural metric on $D$ induced by the Killing form. By the
Cartan decomposition ${\mathfrak{g}}_0={\mathfrak{k}}\oplus
{\mathfrak{p}}_0$ we know that the Killing form $\kappa$ on
${\mathfrak{g}}_0$ is positive
definite on $\mathfrak{p}_0$ and is negative definite on $\mathfrak{k}/%
\mathfrak{v}$. By the identification \eqref{killing10}, for real
tangent
vectors $X,Y\in T_O^{\mathbb{R}} D$, if $X=X_1+X_2$ and $Y=Y_1+Y_2$ where $%
X_1,Y_1\in \mathfrak{k}/\mathfrak{v}$ and $X_2,Y_2\in
\mathfrak{p}_0$, we let
\begin{eqnarray}  \label{killing20}
\tilde\kappa(X,Y)=-\kappa(X_1,Y_1)+\kappa(X_2,Y_2).
\end{eqnarray}
Then $\tilde\kappa$ is a positive definite symmetric bilinear form on $T_O^{%
\mathbb{R}} D$. Now we use left translation of elements in
$G_{\mathbb{R}}$ to move this metric to the real tangent space of
every point in $D$ and we obtain a Riemannian metric on $D$. This is
the Hodge metric defined by Griffiths and Schmid in \cite{grsch1}.

\subsection{Global Torelli Theorem and Applications}

Now we describe the global Torelli theorem.

\begin{theorem}
\label{main} Let $(M,L)$ be a polarized CY manifold of dimension $n$
and let ${\mathcal{T}}$ be its Teichm\"uller space. Let $D$ be the
classifying space of the variation of Hodge structures according to
the middle cohomology of $M$. Let $\p:\T\to D$ be the period map
which maps each point $q\in\T$ to the Hodge decomposition of the
middle dimensional primitive cohomology of $M_q$ which is a point in
$D$. Then the period map ${\Phi}:{\mathcal{T}}\to D$ is injective.
\end{theorem}

Let us describe the main idea of proving the global Torelli theorem.
See \cite{lsty1} for details. In fact we have proved a stronger
result. For any distinct points $p,q\in\T$, in \cite{lsty1} we
showed that the lines $F^n(M_p)$ and $F^n(M_q)$ do not coincide.
This means the first Hodge bundle already determines polarized
marked Calabi-Yau structures.

The first main component in the proof of Theorem \ref{main} is
Theorem \ref{invarkaform}, namely the \ka forms $\omega$ of the
polarized CY metrics are invariant. From this we know that all the
complex structures corresponding to all points in $\T$ are tamed by
$\omega$. If we fix a base point $0\in\T$ then for any point
$q\in\T$ such that $q\ne 0$, the complex structure on $M_q$ is
obtained by deforming the complex structure on $M_0$ via a unique
Beltrami differential $\phi(q)\in A^{0,1}\lb M_0,T_{M_0}^{1,0}\rb$.
Thus we obtained the assigning map
\[
\rho:\T\to A^{0,1}\lb M_0,T_{M_0}^{1,0}\rb
\]
by letting $\rho(q)=\phi(q)$. The assigning map $\rho$ is
holomorphic. In fact, in \cite{lsty4} we proved that the assigning
map $\rho$ is a holomorphic embedding.

According to the work of Todorov \cite{tod1} and our work
\cite{lsty1} we know that $\Omega_q=\sum_{k=0}^n \frac{1}{k!}\lb
\wedge^k \phi(q)\lrcorner\Omega_0\rb$ is a smooth $(n,0)$-form on
$M_q$ where $\Omega_0$ is a properly normalized holomorphic
$(n,0)$-form on $M_0$. It follows from Theorem \ref{flatcoord} and
Lemma \ref{constructn0} that $\Omega_q$ is a holomorphic
$(n,0)$-form on $M_q$. Since the cohomology classes $[\Omega_0]$ and
$[\Omega_q]$ are generators of the Hodge lines $F^n(M_p)$ and
$F^n(M_q)$, it is enough to show that these two classes are not
proportional. Now we look at the Calabi-Yau equation
\[
c_n\Omega_p\wedge\bar\Omega_p= \omega_p^n
\]
where $p\in\T$ is any point, $\Omega_p$ is a properly normalized
nowhere vanishing holomorphic $(n,0)$-form on $M_p$, $\omega_p$ is
the \ka form of the polarized CY metric on $M_p$ and
$c_n=(-1)^{\frac{n(n-1)}{2}}\lb\frac{\sqrt{-1}}{2}\rb^n$. It follows
from Theorem \ref{invarkaform} and the Calabi-Yau equation that if
$[\Omega_q]=c[\Omega_0]$, then $c=1$ and $\phi(q)=0$ which means
that $M_0$ and $M_q$ are isomorphic. This contradicts the assumption
that $q\ne 0$ as points in $\T$ and the global Torelli theorem
follows.

By using the same method we proved the global Torelli theorem of the
\tei space of polarized Hyper-\ka manifolds:
\begin{theorem}
Let $(M,L)$ be a polarized Hyper-\ka manifold and let $\T$ be its
\tei space. Let $D$ be the classifying space of variation of
polarized weight $2$ Hodge structures according to the data of
$(M,L)$. Then the period map $\p:\T\to D$ which maps each point
$q\in\T$ to the Hodge decomposition of the second primitive
cohomology of $M_q$ is injective.
\end{theorem}

In \cite{lsty4} we gave another proof of the global Torelli theorem
of the \tei space of polarized Hyper-\ka manifolds by directly
showing that the cohomology expansion of the canonical $(2,0)$-forms
has no quantum correction. This implies that the Hodge completion of
the \tei space is biholomorphic to the classifying space via the
Harish-Chandra realization.

Another important property of the \tei space is the existence of
holomorphic flat connections. We proved the following theorem in
\cite{lsty4}.

\begin{theorem}
There exists affine structures on the \tei space $\T$ of polarized
CY manifolds. The affine structures are given by global holomorphic
flat connections on $\T$.
\end{theorem}

It is not difficult by using elementary Lie algebra arguments to
establish the affine structure on the complement of the Schubert
cycle which is the nilpotent orbit containing the base point. In our
case we need to prove that the image of the \tei space under the
period map do not intersect certain component of the Schubert cycle
which is codimension one. We call this the partial global
transversality. The problem when partial global transversality holds
is very important in the study of global behavior of the period map.
The Griffiths transversality is too weak to deal with such global
problems of the period map. In \cite{lsty4} we obtained the partial
global transversality by using Yau's solution of the Calabi
conjecture.

As a corollary we proved the holomorphic embedding theorem of the
\tei space of polarized CY manifolds in \cite{lsty4}:
\begin{theorem}
The \tei space of polarized CY manifolds can be holomorphically
embedded into the Euclidean space of same dimension. Furthermore,
its Hodge completion is a domain of holomorphy and there exists a
unique \ke metric on the completion.
\end{theorem}

The methods that we used in \cite{lsty1} can be used to prove the
global Torelli theorem for a large class of manifolds of general
type. Furthermore, these methods can be used to treat the invariance
of the plurigenera for \ka manifolds where the projective case was
proved by Siu.


\end{document}